\documentclass[a4paper,12pt]{article}

\usepackage{mathrsfs}
\usepackage{graphicx}
\usepackage{amsmath}
\usepackage{amsfonts}
\usepackage{amssymb}
\usepackage{amsthm}
\usepackage{textcomp}
\usepackage{color}
\usepackage{geometry}
\usepackage[normalem]{ulem}
\usepackage{placeins}
\usepackage{bm}
\usepackage{wasysym}
\usepackage{titlesec}
\usepackage{txfonts}
\usepackage{graphicx} 
\usepackage{subfigure}

\newtheorem{lemma}{Lemma}[section]

\newtheorem{definition}{Definition}[section]

\renewcommand{\theequation}{\arabic{section}.\arabic{equation}}
\renewcommand{\thetheorem}{\arabic{section}.\arabic{theorem}}
\renewcommand{\thelemma}{\arabic{section}.\arabic{lemma}}
\renewcommand{\theproposition}{\arabic{section}.\arabic{proposition}}

\renewcommand{\thealgorithm}{\arabic{section}.\arabic{algorithm}}

\newcommand\dd{~{\rm d}}

\newcommand\R{\mathbb{R}}

\newcommand\Z{\mathbb{Z}}

\newcommand\RL{\mathcal{R}}

\newcommand\tu{\tilde{u}}

\newcommand\tl{\tilde{\lambda}}
\newcommand\Dn{\widetilde{D}}

\newcommand\intRd{\fint}
\newcommand\Ec{E_{\rm c}}
\newcommand\Nc{N_{\rm c}}
\newcommand\Dos{\mathcal{D}}
\newcommand\Ef{E_{\rm f}}

\title{Plane Wave Methods for Quantum Eigenvalue Problems of Incommensurate Systems}

\author{
	Yuzhi Zhou
	\footnote{{\tt zhou\_yuzhi@iapcm.ac.cn}.
	Institute of Applied Physics and Computational Mathematics, Bejing 100088, China.
	This work was supported by the National Science Foundation of China under grant 91730302.
	},
	~ Huajie Chen
	\footnote{{\tt chen.huajie@bnu.edu.cn}.
	School of Mathematical Sciences, Beijing Normal University, Beijing 100875, China.
	This work was supported by Thousand Talents Program for Young Professionals, and the Fundamental Research Funds for the Central Universities of China under grant 2017EYT22.
	},
	~ and Aihui Zhou
	\footnote{{\tt azhou@lsec.cc.ac.cn}.
	LSEC, Institute of Computational Mathematics and Scientific/Engineering Computing,
	Academy of Mathematics and Systems Science, Chinese Academy of Sciences,  Beijing 100190, China;
	School of Mathematical Sciences, University of Chinese Academy of Sciences, Beijing 100049, China.
	This work was supported by the National Science Foundation of China under grant 91730302 and 11671389, and the Key Research Program of Frontier Sciences of the Chinese Academy of Sciences under grant QYZDJ-SSW-SYS010.}
}
\date{}

\begin{document}
\maketitle

\begin{abstract}
	We propose a novel numerical algorithm for computing the electronic structure related eigenvalue problem of incommensurate systems. 
	Unlike the conventional practice that approximates the system by a large commensurate supercell, our algorithm directly discretizes the eigenvalue problem under the framework of a plane wave method. 
	The emerging ergodicity and the interpretation from higher dimensions give rise to many unique features compared to what we have been familiar with in the periodic system. 
	The numerical results of 1D and 2D quantum eigenvalue problems are presented to show the reliability and efficiency of our scheme. 
	Furthermore, the extension of our algorithm to full Kohn-Sham density functional theory calculations are discussed.
\end{abstract}

\section{Introduction}
\label{sec:intro}
\setcounter{equation}{0}

Recently, there has been growing research interest on the 2D incommensurate layered crystal structures, due to the realization of the heterostructures of 2D materials in experiments and their unique physical properties \cite{britnell13,geim13,Liu16,Novo16,terrones14}. The absence of periodicity presents a fundamental challenge to compute the electronic structure of the incommensurate systems. 
The conventional method to study incommensurate systems is to approximate the systems by a commensurate supercell and then apply Bloch's theory \cite{ebnonnasir14,koda16,komsa13,loh15,terrones14}. 
However, the commensurate supercells are generally very large, thus the electronic structure calculations are computationally expensive or most likely, infeasible.  In addition, the justification of the approximation requires more rigorous studies. 
Recently, there are newly developed approaches without resort to commensurate supercell approximations \cite{cances17,massatt18,massatt17}. These methods are based on the tight-binding models, and the discrete feature of the model cannot be directly generalized to continuous electronic structure models.

The purpose of this paper is to develop a plane wave based numerical framework for solving the (continuous) eigenvalue problems of incommensurate systems. 
The advantage of plane waves is that they form a convenient and efficient representation of the potentials and solutions in each periodic layer.
Even though the whole incommensurate system lacks periodicity, we are still able to discretize the eigenvalue problem with plane waves by exploiting the emerging ergodicity from incommensurate structures. 
Our plane wave discretizations of the incommensurate eigenvalue problem can be interpreted with a periodic higher dimensional form, which help us better format the solution and compute related quantities in full density functional theory (DFT) calculations.
We believe that this numerical framework could lay the foundation for general electronic structure calculations of incommensurate systems in the near future.

The rest of the paper is organized as follow. 
In Section \ref{sec:incommensurate}, we will brifely describe the incommensurate layered systems and the corresponding quantum eigenvalue problem. 
In Section \ref{sec:planewave}, we will first introduce the plane wave discretizations of the eigenvalue problem, then discuss the ergodicity nature and define the density of states, following that the interpretation of the incommensurate eigenvalue problem in higher dimensions re-examine the problem from a different perspective. 
In Section \ref{sec:numerics}, we report some numerical experiments for 1D and 2D quantum eigenvalue problems of some incommensurate systems to show the efficiency of our framework. 
In Section \ref{sec:further}, we further discuss the extension of the current framework to full Kohn-Sham DFT calculations. 
Finally, we provide some concluding remarks in Section \ref{sec:conclusion}.

\section{Incommensurate systems}
\label{sec:incommensurate}
\setcounter{equation}{0}

We consider two $d$-dimensional $(d=1,2)$ periodic systems that are stacked in parallel along the $(d+1)$th dimension. To simplify the presentations, we will neglect the $(d+1)$th dimension and the distance between the two layers. This coordinate is not essential in studying the incommensurate systems and can be easily incorporated into our frameworks. The generalization to incommensurate systems with more than two layers is also straightforward, though with increasing computational cost and complexity.

Each of the $d$-dimensional periodic system can be described by a Bravais lattice
\begin{eqnarray*}
	\RL_j = \big\{ A_j n~:~ n\in\Z^d \big \}, \qquad j=1,2,
\end{eqnarray*}
where $A_j\in\R^{d\times d}$ is invertible. The unit cell for the $j$-th layer is
\begin{eqnarray*}
	\Gamma_j = \big\{ A_j \alpha ~:~ \alpha\in[0,1)^d \big\} , \qquad j=1,2 .
\end{eqnarray*}
The associated reciprocal lattice and reciprocal unit cell are given by
\begin{eqnarray*}
	\RL^*_j &=& \big\{ 2\pi A_j^{-{\rm T}} n~:~ n\in\Z^d \big \} ,
	\\[1ex]
	\Gamma^*_j &=& \big\{ 2\pi A_j^{-{\rm T}}  \alpha ~:~ \alpha\in[0,1)^d \big\}
\end{eqnarray*}
respectively, for $j=1,2$.

Although each individual lattice $\RL_j$ is periodic in the sense that
\begin{eqnarray*}
	\RL_j = A_j n + \RL_j \qquad \forall~n\in \Z^d,~~j=1,2 ,
\end{eqnarray*}
the joined system $\RL_1\cup\RL_2$ need not to be periodic.
We consider the {\it incommensurate} system defined as follows.

\begin{definition}[Incommensurateness]
	\label{incomm}
	Two lattices $\RL_1$ and $\RL_2$ are incommensurate if
	\begin{eqnarray}
	\RL_1^*\cup\RL_2^* + \tau = \RL_1^*\cup\RL_2^* \quad \Leftrightarrow \quad \tau=\pmb{0}\in\R^d .
	\end{eqnarray}
	Otherwise, the lattices $\RL_1$ and $\RL_2$ are commensurate.
\end{definition}

We consider the following Schr\"{o}dinger-type eigenvalue problem for an incommensurate system:
Find $(\lambda,u)$ such that
\begin{eqnarray}\label{eq:eigen}
\Big( -\frac{1}{2}\Delta + V_1({\bf r}) + V_2({\bf r}) \Big) u({\bf r}) = \lambda u({\bf r}) \qquad{\rm for}~{\bf r}\in\R^d ,
\end{eqnarray}
where $V_j:\R^d\rightarrow\R$ is $\RL_j$-periodic
\begin{eqnarray}\label{V_periodic}
V_j({\bf r}+\tau)=V_j({\bf r}) \qquad \forall~\tau\in\RL_j \quad {\rm for}~ j=1,2.
\end{eqnarray}
We assume throughout this paper that $V_j$ are smooth $(j=1,2)$ and the lattices $\RL_1$ and $\RL_2$ are incommensurate.
Note that \eqref{V_periodic} implies that the potentials $V_j$ can be written as Fourier series
\begin{eqnarray}\label{V_series}
V_j({\bf r}) = \sum_{m\in\Z^d} V_{jm} e^{iG_{jm}\cdot {\bf r}}
\qquad{\rm for} ~j=1,2,
\end{eqnarray}
where $G_{jm}=A_j^{\rm -T}m\in\RL_j^*$ are wavevectors in the reciprocal lattice, and
\begin{eqnarray*}
	V_{jm} = \frac{1}{|\Gamma_j|}\int_{\Gamma_j}V_j({\bf r})e^{-iG_{jm}\cdot{\bf r}}\dd {\bf r}.
\end{eqnarray*}

Solving \eqref{eq:eigen} is the central task of many electronic structure related quantum models, for example,
Gross-Pitaevskii equations \cite{gross61}, Kohn-Sham equations \cite{martin04},
and Hartree-Fock equations \cite{helgaker00}. 
For periodic systems, Bloch's theorem decomposes and diagonalizes the eigenvalue problem by the wavevectors within the first Brillouin zone. 
For incommensurate layered systems, there is no periodicity overall. 
However, the periodicity within each layer and the incommensurate nature between layers impose unique features on the problem and solution,
which will help us design numerical methods for the problem.

\section{Plane wave discretizations}
\label{sec:planewave}
\setcounter{equation}{0}

\subsection{General expressions}

As mentioned earlier, the plane wave representation is especially appropriate for description of the periodic potentials.
We shall generate the representation to incommensurate layered systems.

We first introduce the following notations.
Denote the average spacial integral by 
\begin{eqnarray*}
\intRd:=\lim_{R\rightarrow\infty} \frac{1}{|B_R|}\int_{B_R} ~ ,
\end{eqnarray*}
where $B_R\subset\R^d$ is the ball centred at origin with radii $R$.
We immediately have the following orthonormal condition:
\begin{eqnarray}\label{ortho}
\intRd e^{-i{\bf k}{\bf r}} e^{i{\bf k}'{\bf r}} \dd{\bf r} = \delta_{{\bf k}{\bf k}'}
\qquad\forall~{\bf k},{\bf k}'\in\R^d.
\end{eqnarray}
For appropriate function $u:\R^d\rightarrow\R$, we define the following {\it averaged} Fourier transform:
\begin{eqnarray}\label{fourier}
\hat{u}({\bf k}):=\intRd u({\bf r})e^{-i{\bf k}\cdot {\bf r}} \dd {\bf r} \qquad {\rm for}~{\bf k} \in\R^d.
\end{eqnarray}
The above definition depends on the existence of limit on the right-hand side, and
we shall provide a rigorous functional space for this in a forthcoming math paper \cite{chen19}.

We seek eigenfunctions of \eqref{eq:eigen} by performing the above transform
\begin{eqnarray}\label{transform}
\intRd e^{-i{\bf k}\cdot{\bf r}} \Big( -\frac{1}{2}\Delta + V_1({\bf r}) + V_2({\bf r}) - \lambda  \Big) u({\bf r})  \dd {\bf r} = 0
\qquad{\rm for}~{\bf k}\in\R^d.
\end{eqnarray}
This together with \eqref{V_series} and \eqref{ortho} implies
\begin{eqnarray}\label{eq:eigen-k}
\frac{1}{2}|{\bf k}|^2 \hat{u}({\bf k}) + \sum_{m\in\Z^d}V_{1m}\hat{u}({\bf k}-G_{1m})
+ \sum_{n\in\Z^d}V_{2n}\hat{u}({\bf k}-G_{2n})  = \lambda\hat{u}({\bf k})
\qquad {\rm for}~{\bf k}\in\R^d .
\end{eqnarray}

It can been seen from \eqref{eq:eigen-k} that these equations couple the wavevectors ${\bf k}$ and ${\bf k'}$
only if they differ by a sum of \emph{two} lattice vectors $G_{1m}\in\RL_1^*$ and $G_{2n}\in\RL_2^*$:
\begin{eqnarray}\label{couple_k}
{\bf k}-{\bf k}' = G_{1m}+G_{2n} = 2\pi\big( A_1^{\rm -T}m + A_2^{\rm -T}n \big)
\qquad ~{\rm for~some}~ m,n\in\Z^d .
\end{eqnarray}
This is very much similar to the case in periodic systems,  except that ${\bf k}$ and ${\bf k'}$ are only coupled by a \emph{single} lattice vector in the periodic systems. 
With the coupling relations between wavevectors, the eigenvalue problem \eqref{eq:eigen}
for any given ${\bf k}\in\R^d$ can be written as an (infinite) matrix equation
\begin{eqnarray}\label{eq:eigen-matrix}
\sum_{m',n'\in\Z^d} H_{mn,m'n'}({\bf k}) \hat{u}({\bf k}+G_{1m'}+G_{2n'})
= \lambda({\bf k}) \hat{u}({\bf k}+G_{1m}+G_{2n})
\qquad m,n\in\Z^d,
\end{eqnarray}
where
\begin{eqnarray}\label{H_elements}
H_{mn,m'n'}({\bf k}) = \frac{1}{2} |{\bf k}+G_{1m}+G_{2n}|^2 \delta_{mm'}\delta_{nn'} +
V_{1(m-m')}\delta_{nn'} + V_{2(n-n')}\delta_{mm'} .
\end{eqnarray}
While the periodic problem can sample the ${\bf k}$-points in the first Brillouin zone with Bloch's theory,
the choice of ${\bf k}$-point in \eqref{eq:eigen-matrix} depends on how the wavevectors $\big\{{\bf k}+G_{1m}+G_{2n}\big\}$ distribute in reciprocal space.
This will be discussed in the following subsection.

Note that the above formulas only hold for incommensurate systems.
Since for commensurate systems, there exist ${\bf 0}\neq G_{1m}\in\RL_1^*$ and ${\bf 0}\neq G_{2n}\in\RL_2^*$ such that $G_{1m}+G_{2n}={\bf 0}$,
in which case the expression of \eqref{eq:eigen-matrix} becomes redundant and does not correspond to the eigenvalue problem of a periodic system.

\subsection{Ergodicity and density of states}
\label{sec:ergodicityDOS}

This essential difference between periodic and incommensurate systems is the so-called {\it ergodicity}.
This term is originally from statistical mathematics and thermodynamics, which describes the equiprobable access to all states in the phase space.
The ergodicity is a direct consequence from the incommensurateness defined in \eqref{incomm}, and is the root of unique features of the incommensurate eigenvalue problem.
It can be stated in the mathematical language as those in \cite[Proposition 2.4]{cances17}, \cite[Theorem 2.1]{massatt17} and \cite{dingzhou09}.

\begin{lemma}[Ergodicity]\label{lemma:ergodic}
	If $\RL_1$ and $\RL_2$ are incommensurate lattices, then the set
	$\big\{ A_1^{\rm -T}m+A_2^{\rm -T}n \big\}_{m,n\in\Z^d}$
	is dense and uniformly distributed in $\R^d$.
\end{lemma}

The ergodic nature is reflected in both real and reciprocal spaces. We first discuss the ergodicity in a reciprocal space. 
The ergodicity in a real space will be discussed in the next subsection.

\begin{figure}[hbt!]
	\centering
	\includegraphics[width=15cm]{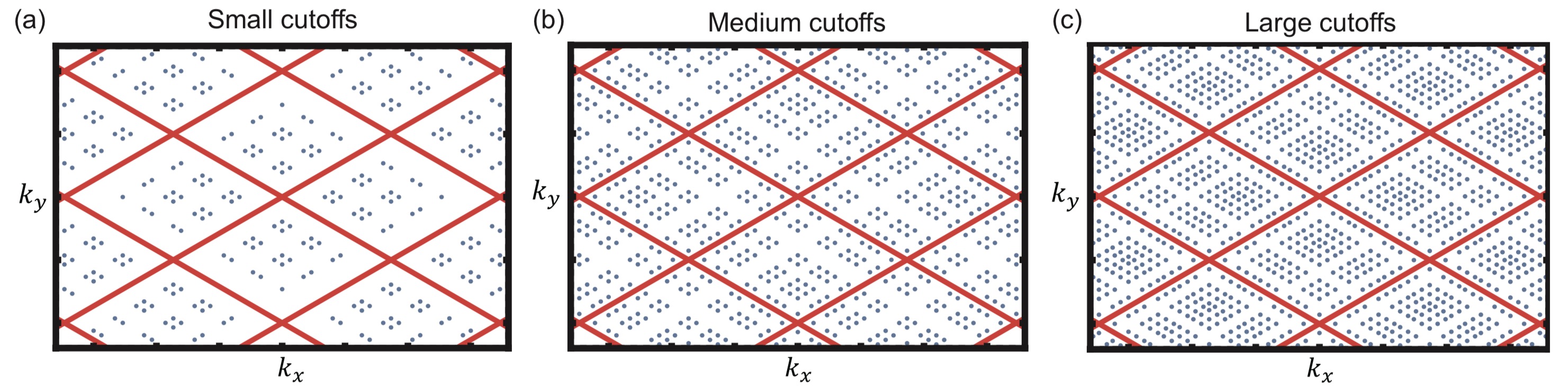}
	\caption{The illustration of ergodicity in the reciprocal space. The 2D incommensurate system consists of two hexagonal lattices with reciprocal lattice constants $1$ and $\sqrt{3}$, respectively. 
	The red grid lines divide the reciprocal space into the hexagonal Brillouin zones whose lattice constant is $1$. 
	The purpose of these pictures is to better visualize the distribution of wavevectors $\{{\bf k}+G_{1m}+G_{2n}\}$.}
	\label{fig:ErgoInRec}
\end{figure}

Lemma \ref{lemma:ergodic} implies that for any single ${\bf k}\in\R^d$, its coupled wavevectors $\big\{{\bf k}+G_{1m}+G_{2n}\big\}$ could densely and uniformly spread out the reciprocal space, as $m,n\rightarrow\infty$. 
This is visualized in Fig.~\ref{fig:ErgoInRec}, where the reciprocal space of the incommensurate hexagonal 2D lattice is sampled by the wavevector set ${\bf k}+G_{1m}+G_{2n}$ generated from $\Gamma$ point (${\bf k}=\pmb{0}$) and different cutoffs of $m$ and $n$. 
It can be seen that as the cutoffs of $m,n$ increase, the wavevector set gradually becomes dense and uniform in the reciprocal space.

This fact has two implications. First, the dense wavevector set implies that the spectrum of \eqref{eq:eigen-matrix} could be continuous.
Second, the full spectrum structure of \eqref{eq:eigen}, in principle, can be adequately reconstructed by solving \eqref{eq:eigen-matrix} with single ${\bf k}$-point and sufficiently large cutoffs. 
This is significantly different from the case of periodic systems.
For periodic problems, the wavevector set generated by a single ${\bf k}$-point only contains discrete points in the reciprocal space, and only point spectrum can be obtained.

In practical simulations, we need to restrict \eqref{eq:eigen-matrix} to a finite dimensional subspace with some energy cutoff $\Ec$. 
More precisely, we require $|G_{1m}|^2+|G_{2n}|^2\leq 2\Ec$ and obtain a finite set of plane wave vectors, with $\Nc$ the number of wavevectors in the set.
We shall first restrict ourself to the calculations with one single ${\bf k}$-point.
Then we can obtain a discrete set of eigenvalues $\lambda_j({\bf k}),~j=1,\cdots,\Nc$ by solving the matrix eigenvalue problem (which is a truncated form of \eqref{eq:eigen-matrix})
\begin{multline*}
	\qquad
	\sum_{{m', n'\in\Z^d}\atop{|G_{1m'}|^2+|G_{2n'}|^2\leq 2\Ec}} H_{mn,m'n'}({\bf k}) \hat{u}_j({\bf k}+G_{1m'}+G_{2n'})
	= \lambda_j({\bf k}) \hat{u}_j({\bf k}+G_{1m}+G_{2n})
	\\
	{\rm for}~ m,n\in\Z^d ~~{\rm and}~~ |G_{1m}|^2+|G_{2n}|^2\leq 2\Ec ,
	\qquad
\end{multline*}
where the matrix elements $H_{mn,m'n'}({\bf k})$ are given by \eqref{H_elements}.
The approximation of the $j$-th eigenfunction can be written as
\begin{eqnarray}\label{eigenfunction}
u_{j,{\bf k}}({\bf r}) = \sum_{{m, n\in\Z^d}\atop{|G_{1m}|^2+|G_{2n}|^2\leq 2\Ec}} \hat{u}_{j}({\bf k}+G_{1m}+G_{2n}) e^{i({\bf k}+G_{1m}+G_{2n})\cdot{\bf r}} .
\end{eqnarray}
Using \eqref{ortho}, we see that $u_{j,{\bf k}}({\bf r})$ is automatically normalized with respect to the norm
$\|u\|:=\left(\intRd u^2\right)^{1/2}$.

To represent the (continuous) spectrum sturcture, we use the following definition of {\it density of states} (DoS):
\begin{eqnarray}\label{def:dos}
\Dos(\epsilon,{\bf k},\Ec) := \frac{1}{\sqrt{\Ec}} \sum_{j=1}^{\Nc}\delta\big(\epsilon-\lambda_j({\bf k})\big) .
\end{eqnarray}
Here, $\frac{1}{\sqrt{\Ec}}$ is a normalization prefactor such that $\Dos$ can represent the averaged DoS,  
and we refer to Section \ref{sec:further} for detailed discussions of the choice of the prefactor.
This definition will be slightly improved in Section \ref{sec:further} to enclose more physical meaning.

We expect that the DoS does not depend on the choice of ${\bf k}$ and can converge to the real DoS of \eqref{eq:eigen}
as $\Ec$ goes to infinity:
\begin{eqnarray}
\Dos(\epsilon,{\bf k},\Ec) \xrightarrow{\Ec\rightarrow\infty} \Dos(\epsilon) \qquad\forall~{\bf k}\in\R^d.
\end{eqnarray}
We will support this convergence with some numerical examples in Section \ref{sec:numerics},
and also provide a rigorous proof in a forthcoming mathematical paper \cite{chen19}.

When it comes to more practical total energy calculations, sampling the reciprocal space with multiple ${\bf k}$-points must be conducted.
We will discuss the principles of ${\bf k}$-point sampling to maximally recover the continuous spectrum of \eqref{eq:eigen} in the following subsection.

\subsection{Interpretation in higher dimensions}
\label{sec:highDim}

In this section, we will show that the incommensurate eigenvalue problem can be interpreted in a higher dimension, in which the periodicity can be restored.
We mention that similar idea has been explored for the lattices and diffraction patterns of quasi-crystals (see e.g. \cite{baake17,blionv15,jiang14,Walter09}).

For the incommensurate eigenvalue problem \eqref{eq:eigen}, we construct the following eigenvalue problem in a higher dimensional space $\R^{d}\times\R^d$:
\begin{eqnarray}\label{eq:eigen-H}
\Big( -\frac{1}{2}\Dn + V_1({\bf r}) + V_2({\bf r}') \Big) \tu({\bf r},{\bf r}') = \tl \tu({\bf r},{\bf r}')
\qquad{\rm for}~({\bf r},{\bf r}')\in\R^d\times\R^d ,
\end{eqnarray}
where the differential operator $\Dn$ is definded by
\begin{eqnarray*}
	\Dn\tilde{u}({\bf r},{\bf r}') := \sum_{i=1}^d \big(\partial_{{\bf r}_i} + \partial_{{\bf r}'_i}\big)^2 \tilde{u}({\bf r},{\bf r}')
	\qquad ({\bf r},{\bf r}')\in\R^d\times\R^d.
\end{eqnarray*}
Since $V_1$ and $V_2$ are periodic in $\R^d$ with respect to $\RL_1$ and $\RL_2$ respectively,
the potential $V_1({\bf r}) + V_2({\bf r}')$ is periodic in $\R^d\times\R^d$ with respect to the higher dimension lattice
\begin{eqnarray*}
	\widetilde{\RL} := \RL_1 \times \RL_2
	= \Big\{ \big(A_1^{\rm T}m,A_2^{\rm T}n\big) ~:~ (m,n)\in\Z^d\times\Z^d \Big\} .
\end{eqnarray*}
Therefore, the operator $\widetilde{H} = -\frac{1}{2}\Dn + V_1(x) + V_2(y)$ is translation invariant with respect to the lattice $\widetilde{\RL}$.
Hence \eqref{eq:eigen-H} is a periodic problem in $\R^d\times\R^d$ and we can apply Bloch's theory to it.
Note that the resolvent of operator $\widetilde{H}$ is not compact on the cell
$\Gamma_1\times\Gamma_2\subset\R^d\times\R^d$ (with the periodic boundary condition).
Therefore, the spectrum of \eqref{eq:eigen-H-matrix} is not a discrete set.
As the energy cutoff for plane wave vectors goes to infinity, the spectrum can become continuous.
The constructions of DoS in the previous section can be directly extended to this higher dimensional problem.

Let $\widetilde{\RL}^*$ be the reciprocal lattice of $\widetilde{\RL}$,
and $\widetilde{\Gamma}^*$ be the unit cell of $\widetilde{\RL}^*$.
By applying Bloch's theorem, we can derive that for a given $\widetilde{\bf k} = (\tilde{\bf k}_1,\tilde{\bf k}_2)\in\widetilde{\Gamma}^*$,
the eigenstates of \eqref{eq:eigen-H}  at $\widetilde{\bf k}$, denoted by $\tilde{\lambda}(\tilde{\bf k})$ and $\tilde{u}_{\tilde{\bf k}}$,
can be obtained by solving
\begin{eqnarray}\label{eq:eigen-H-matrix}
\sum_{m',n'\in\Z^d} \widetilde{H}_{mn,m'n'}\big(\widetilde{\bf k}\big) U_{m'n'}(\widetilde{\bf k})
= \tilde{\lambda}(\widetilde{\bf k}) U_{mn}(\widetilde{\bf k})
\qquad m,n\in\Z^d,
\end{eqnarray}
where
\begin{eqnarray}\label{H_tilde_elements}
\widetilde{H}_{mn,m'n'}(\widetilde{\bf k}) 
= \frac{1}{2} \big|\tilde{\bf k}_1+\tilde{\bf k}_2+G_{1m}+G_{2n}\big|^2 \delta_{mm'}\delta_{nn'} +
V_{1(m-m')}\delta_{nn'} + V_{2(n-n')}\delta_{mm'} 
\end{eqnarray}
and the eigenfunction $\tilde{u}_{\tilde{\bf k}}$ can be written as
\begin{eqnarray}\label{eigenfunction-high-dimension}
\tilde{u}_{\tilde{\bf k}}({\bf r},{\bf r}') = \sum_{m,n\in\Z^d}U_{mn}(\tilde{\bf k})
\exp \Big( i (\widetilde{\bf k} + \widetilde{G}_{mn})\cdot ({\bf r},{\bf r}') \Big)
\quad{\rm with}\quad
\widetilde{G}_{mn} = (G_{1m},G_{2n}) \in \widetilde{\RL}^*.
\quad
\end{eqnarray}
We observe that \eqref{eq:eigen-H-matrix} is exactly the same as \eqref{eq:eigen-matrix} by taking ${\bf k}=\tilde{\bf k}_1+\tilde{\bf k}_2$ in \eqref{eq:eigen-matrix},
and hence gives the same spectrum and DoS with a given $\widetilde{\bf k}$. 

Due to the ergodicity in Lemma \ref{lemma:ergodic} and the fact that $\widetilde{H}_{mn,m'n'}(\widetilde{\bf k}) $ depends only on $\tilde{\bf k}_1+\tilde{\bf k}_2$,
we see that the full spectrum can be restored from one single $\widetilde{\bf k}$ with sufficiently large cutoffs of $m$ and $n$.
Alternatively, we can sample ${\bf k}$-points uniformly in the first Brillouin zone of either lattice $\RL_1$ or $\RL_2$ or $\widetilde{\RL}$.
For multiple ${\bf k}$-points case, we need to add an additional prefactor $\frac{1}{N_k}$ in front of the definition of DoS \eqref{def:dos}, 
where $N_k$ is the number of ${\bf k}$-points sampled in the first Brillouin zone.

With the above construction, we can transform the incommensurate problem into a periodic problem in a higher dimension. 
The eigenstates of \eqref{eq:eigen} and \eqref{eq:eigen-H} are identical with the relation $\lambda=\tilde{\lambda}$ and $u({\bf r})=\tilde{u}({\bf r},{\bf r})$.

The higher dimensional (periodic) interpretation not only gives us an alternative expression of the incommensurate eigenvalue problem,
but also facilitates the formulation of full DFT calculations, see Section \ref{sec:further}. 
Moreover, the higher dimensional interpretation can manifest the ergodic nature in real space, which will be discussed in the following.

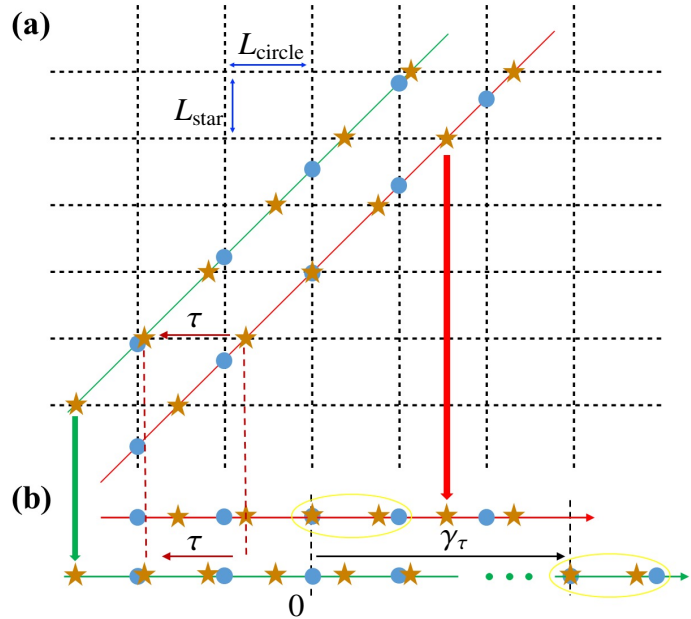
\begin{figure}[hbt!]
	\centering
	\includegraphics[width=10cm]{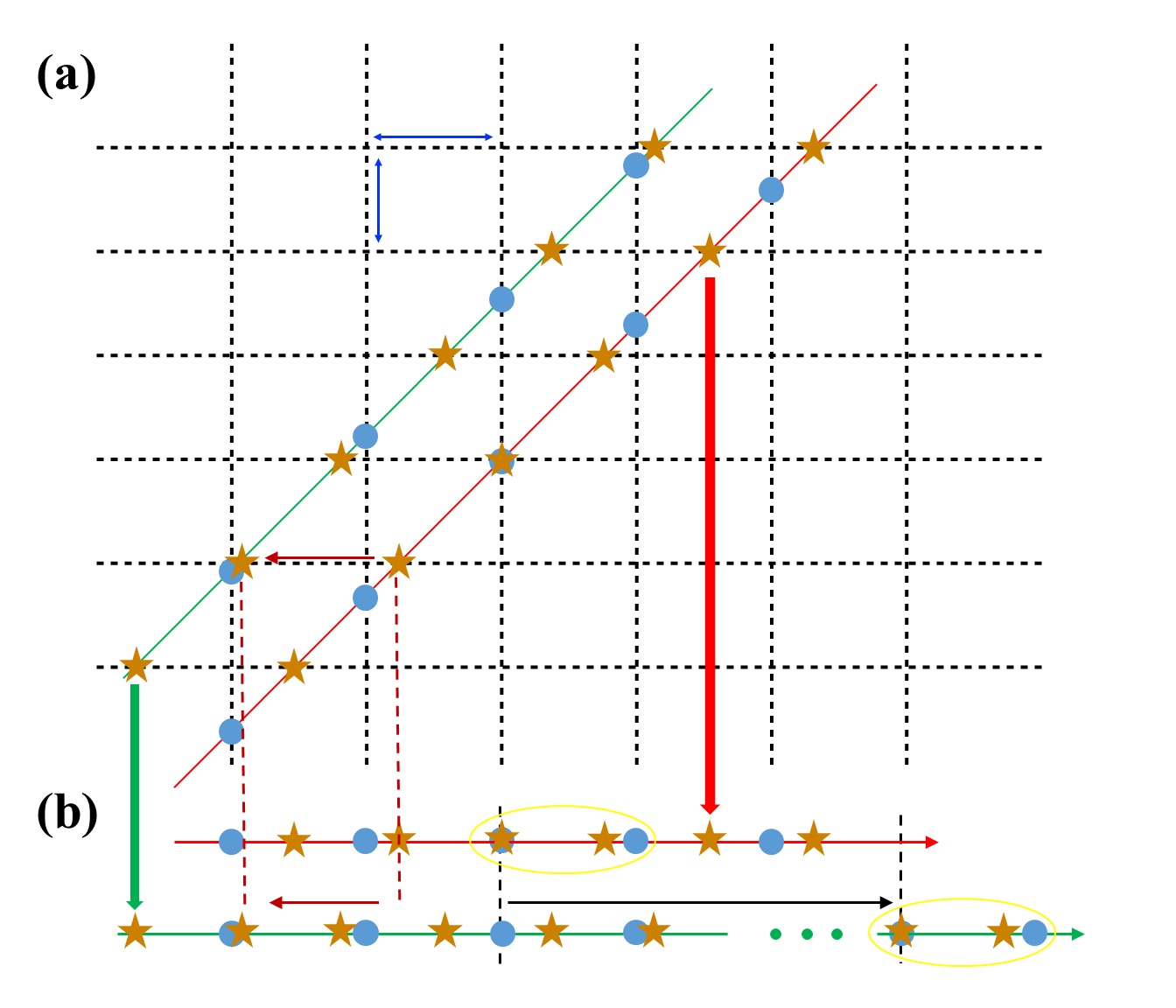}
	\put(-215,190){$L_{\rm star}$}
	\put(-191,215){$L_{\rm circle}$}
	\put(-210,28){$\tau$}
	\put(-210,112){$\tau$}
	\put(-115,30){$\gamma_{\tau}$}
	\put(-172,3){$0$}
	\caption{
		The illustration of ergodicity in the real space. The 1D incommensurate system consists of two atomic chains (ploted with blue circles and red stars) 
		with lattice constants $L_{\rm circle}$ and $L_{\rm star}$.
		(a) The connection between the higher dimensional interpretation and the original incommensurate problem.
		(b) The ergodicity in the real space. After shifting the second atomic chains (red stars) by $\tau$, 
		we can find a vector $\gamma_{\tau}$, such that the local atomic arrangement at $\gamma_{\tau}$ can be almost the same as that at the origin in the original system. }
	\label{fig:ErgoInReal}
\end{figure}

As a matter of fact, the solution of the higher dimensional eigenvalue problem \eqref{eq:eigen-H-matrix} provides more informations than the solution of \eqref{eq:eigen}. 
It contains the solutions for a series of incommensurate systems, which are generated by shiftings one layer with respect to the other.
For example, we can restrict \eqref{eq:eigen-H} on the subset $\big\{({\bf r},{\bf r}'): {\bf r}+\tau = {\bf r}' \big\}$ with a vector $\tau\in\R^d$, and derive the following eigenvalue problem:
\begin{eqnarray}\label{eq:eigen-d}
\Big( -\frac{1}{2}\Delta + V_1({\bf r}) + V_2({\bf r}+\tau) \Big) u_{\tau}({\bf r}) = \tl u_{\tau}({\bf r})
\qquad{\rm for}~{\bf r}\in \R^d 
\end{eqnarray}
with $u_{\tau}({\bf r}) = \tilde{u}({\bf r},{\bf r}+\tau)$.
This is again an incommensurate eigenvalue problem, which is similar to \eqref{eq:eigen}, only that the second lattice is shifted by $\tau$ (if $\tau=0$, then \eqref{eq:eigen-d} is identical to \eqref{eq:eigen}). 
Due to the ergodicity of incommensurate systems, it is nature that the problems with any $\tau\in\R^d$ are {\rm almost} the same with each other, and hence share the same spectrum structure.
More precisely, if one of the layer is shifted by $\tau$, then within arbitrarily required precision, we can find a translation vector $\gamma_{\tau}$ such that at the local atomic configuration at $\gamma_{\tau}$ of the shifted system is the same as that at the origin of the unshifted system. 
This can be visualized for an 1D incommensurate system with two atomic chains in Fig. \ref{fig:ErgoInReal}. 

Note that the above constructions can not be applied to the commensurate systems.
Without ergodicity, shifting one of the lattice will not restore the structure of the original system most of the time, and therefore may change the spectrum.

\section{Numerical simulations}
\label{sec:numerics}
\setcounter{equation}{0}

In this section, we will present the numerical simulations of some linear eigenvalue problems from 1D and 2D incommensurate systems, by using our plane wave methods.
To present the DoS with smooth curves, we use a normalized Gaussian $c_{\rm n}\exp\big(-\sigma(\epsilon-\lambda)^2\big)$ (with $\sigma=5.0$ and $c_{\rm n}$ the normalization constant)
to smear the Dirac function $\delta(\epsilon-\lambda)$ in \eqref{def:dos}.

\begin{figure}[!htb]
	\begin{minipage}[t]{0.5\linewidth}
		\centering
		\includegraphics[width=7cm]{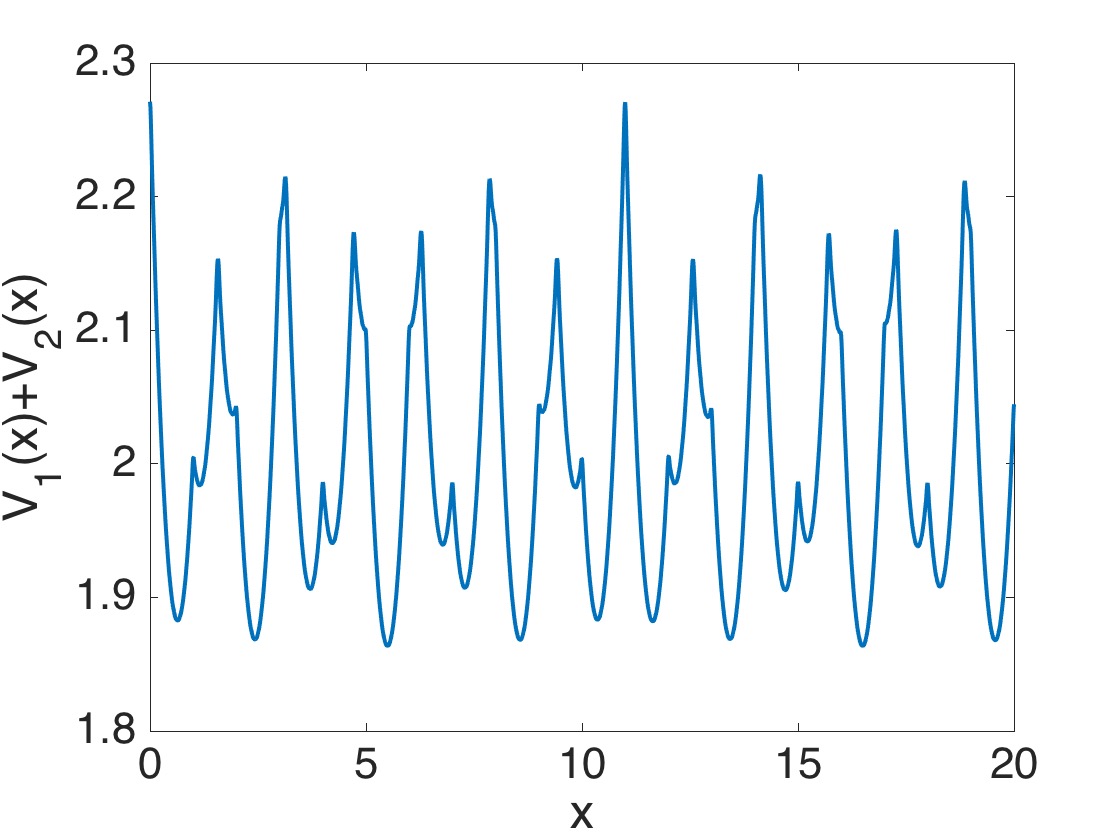}
		\caption{1D incommensurate potential.}
		\label{fig:example1:potential}
	\end{minipage}
	\hspace{0.3cm}
	\begin{minipage}[t]{0.5\linewidth}
		\centering
		\includegraphics[width=7cm]{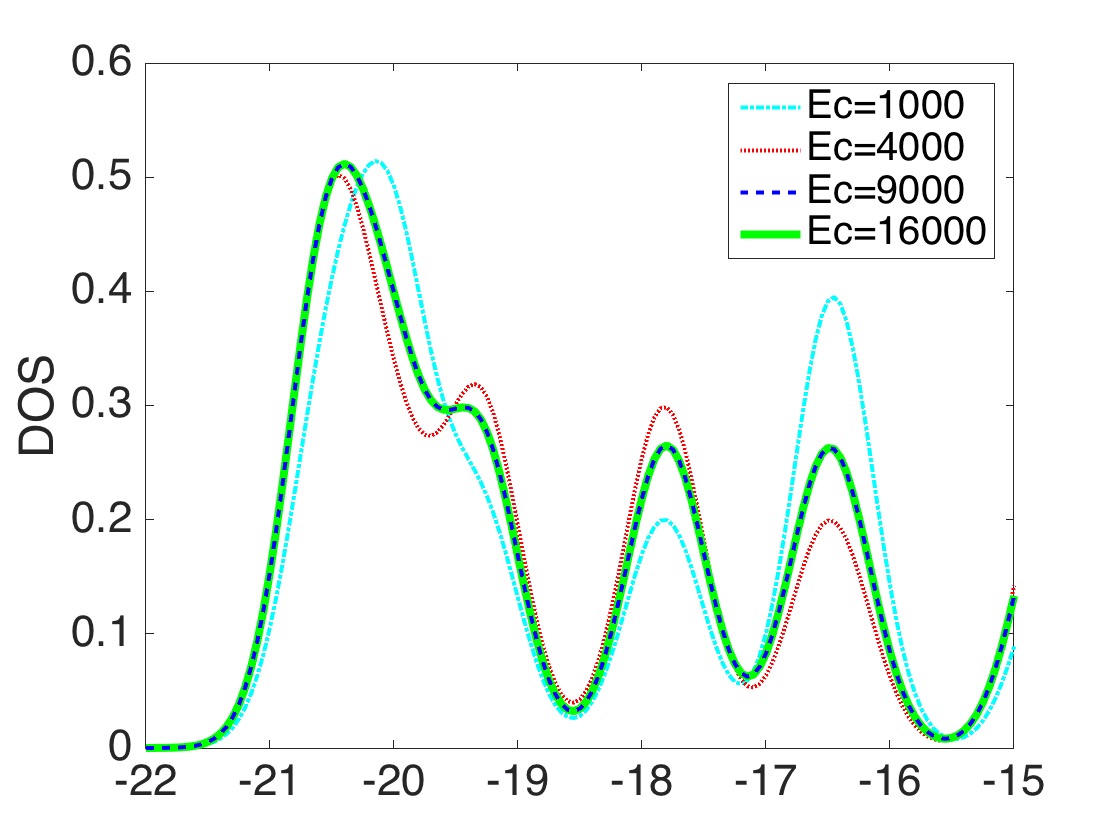}
		\caption{Convergence of $\Dos(\epsilon)$ with respect to the energy cutoff $\Ec$.}
		\label{fig:example1:dos:convergence_Ec}
	\end{minipage}
\end{figure}

{\bf Example 1.}
(one-dimensional chains with incommensurate lattice constants).
Consider the following eigenvalue problem:
\begin{eqnarray}\label{example1}
-u''(x) + \big( V_1(x) + V_2(x) \big)u(x) = \lambda u(x) \qquad x\in\R,
\end{eqnarray}
where $V_1$ and $V_2$ are screened Coulomb potentials with different periodicity
\begin{eqnarray}\label{coulomb}
V_1(x) = Z_1 \sum_{G_{1m}\in \frac{2\pi}{L_1}\Z} \frac{e^{iG_{1m}x}}{|G_{1m}|^2+z}
\quad{\rm and}\quad
V_2(x) = Z_2 \sum_{G_{2n}\in \frac{2\pi}{L_2}\Z} \frac{e^{iG_{2n}x}}{|G_{2n}|^2+z}
\end{eqnarray}
with $L_1=1$, $L_2=\pi/2$, $Z_1=Z_2=1$ and  $z=1$.
The incommensurate potential $V_1(x) + V_2(x)$ is shown in Fig. \ref{fig:example1:potential}.

\begin{figure}[!htb]
	\begin{minipage}[t]{0.5\linewidth}
		\centering
		\includegraphics[width=7cm]{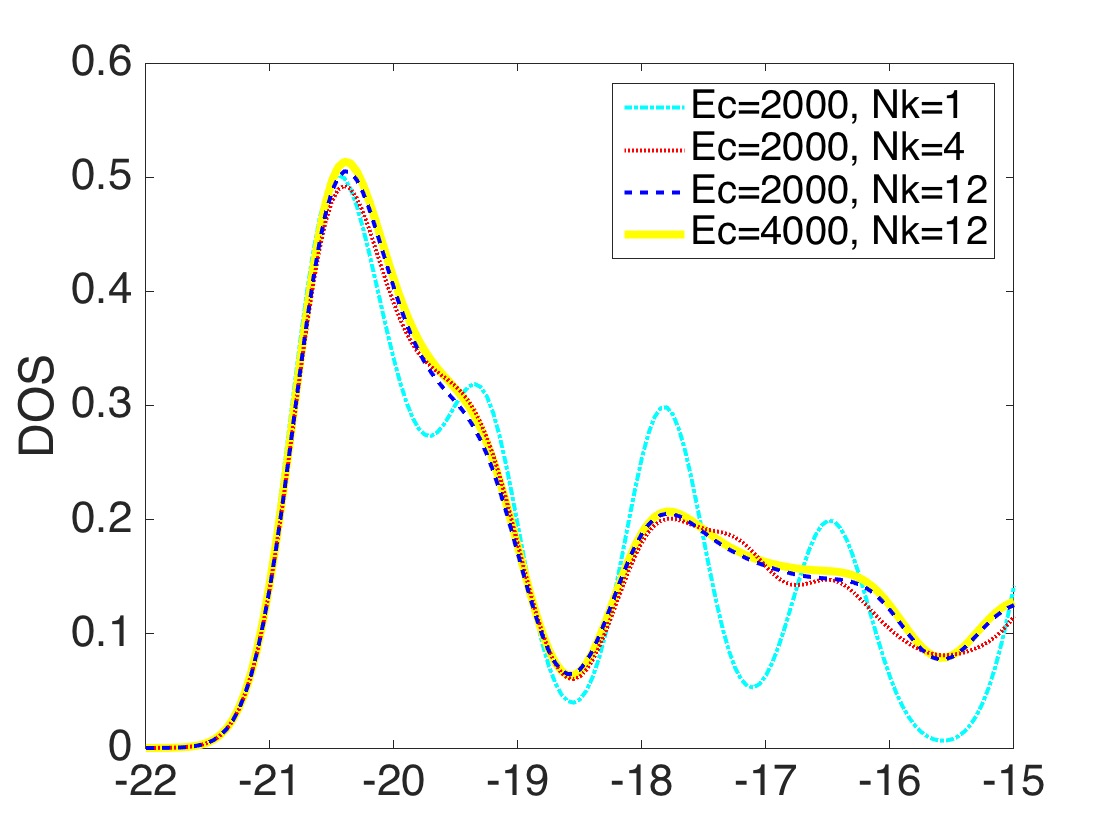}
		\caption{Convergence of $\Dos(\epsilon)$ with respect to ${\bf k}$-sampling.}
		\label{fig:example1:dos:convergence_k}
	\end{minipage}
	\hspace{0.3cm}
	\begin{minipage}[t]{0.5\linewidth}
		\centering
		\includegraphics[width=7cm]{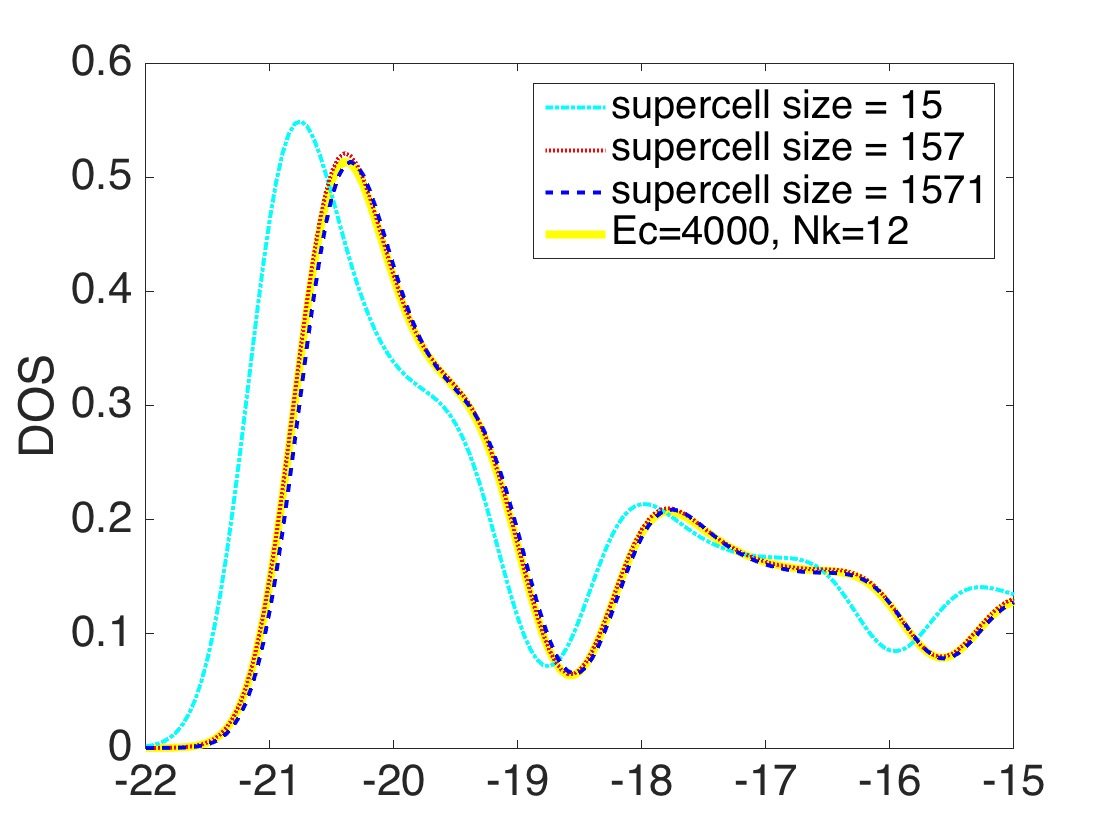}
		\caption{Convergence of DoS from commensurte supercell approximations with
			$L_1=1$ and $L_2=1.5$ (error $\approx 0.07$), $L_2=1.57$ (error $\approx 0.001$), $L_2=1.571$ (error $\approx 0.0002$), respectively.}
		\label{fig:example1:dos:commensurate}
	\end{minipage}
\end{figure}

We use a single $\Gamma$ point (${\bf k}={\bf 0}$) and different energy cutofs $\Ec$ to solve \eqref{example1}.
The DoS are shown in Fig. \ref{fig:example1:dos:convergence_Ec},
from which we observe that the convergence with respect to energy cutoff $\Ec$.
Furthermore, we repeat the simulations with multiple ${\bf k}$-points with a given $\Ec$,
and show the convergence of DoS in Fig. \ref{fig:example1:dos:convergence_k}.
With comparisons, we observe that sampling more ${\bf k}$-points could be more efficient 
to achieve convergence than simply increasing $\Ec$ with a single ${\bf k}$-point.
We point out that the two DoS limits look slightly different due to the smearing width of the Gaussian (used to plot DoS) and the fact that
a single ${\bf k}$-point does not creat enough eigenstates in the high energy window (when $\Ec$ is not large enough, and hence converge slower than that using multiple ${\bf k}$-points).

We also use the standard supercell approximation method to simulate this incommensurat problem.
The DoS obtained from commensurate supercell approximations are presented in Fig. \ref{fig:example1:dos:commensurate}, 
from which we see that a very large supercell must be used in the simulation to achieve similar accuracy.

\begin{figure}[htb!]
	\centering
	\includegraphics[width=4.8cm]{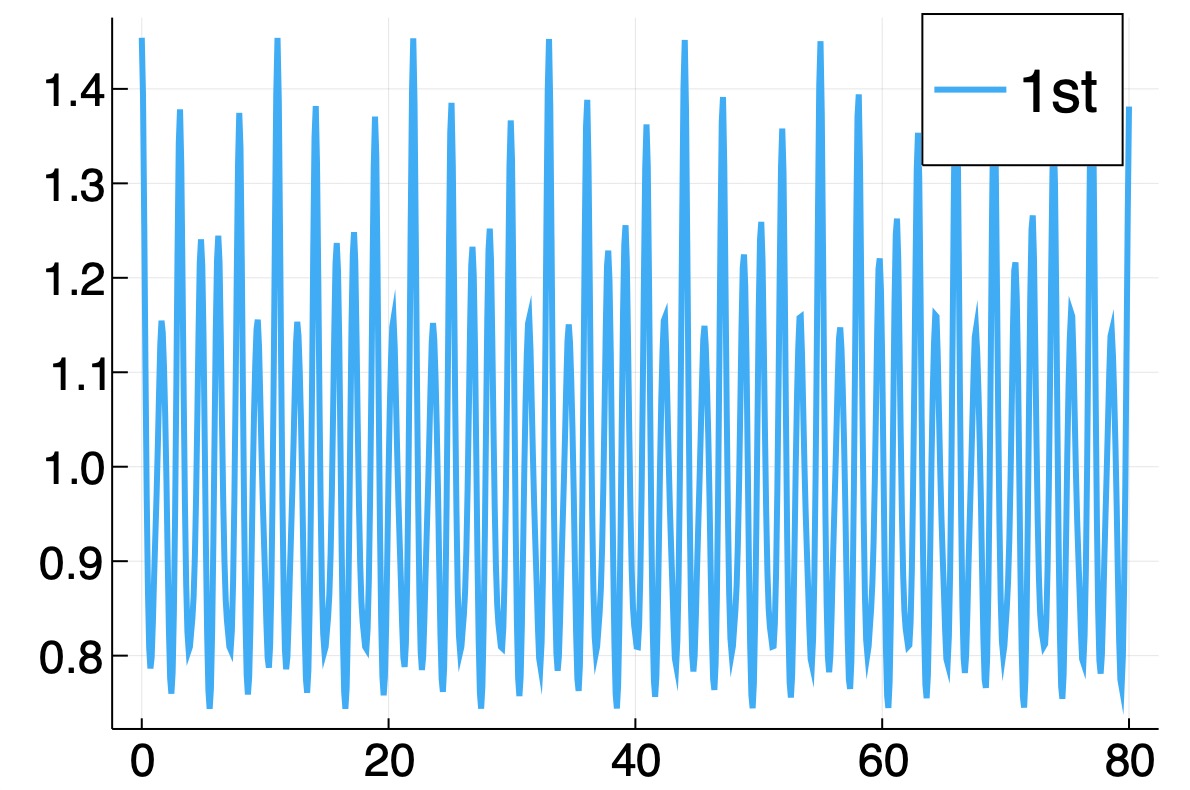}
	\includegraphics[width=4.8cm]{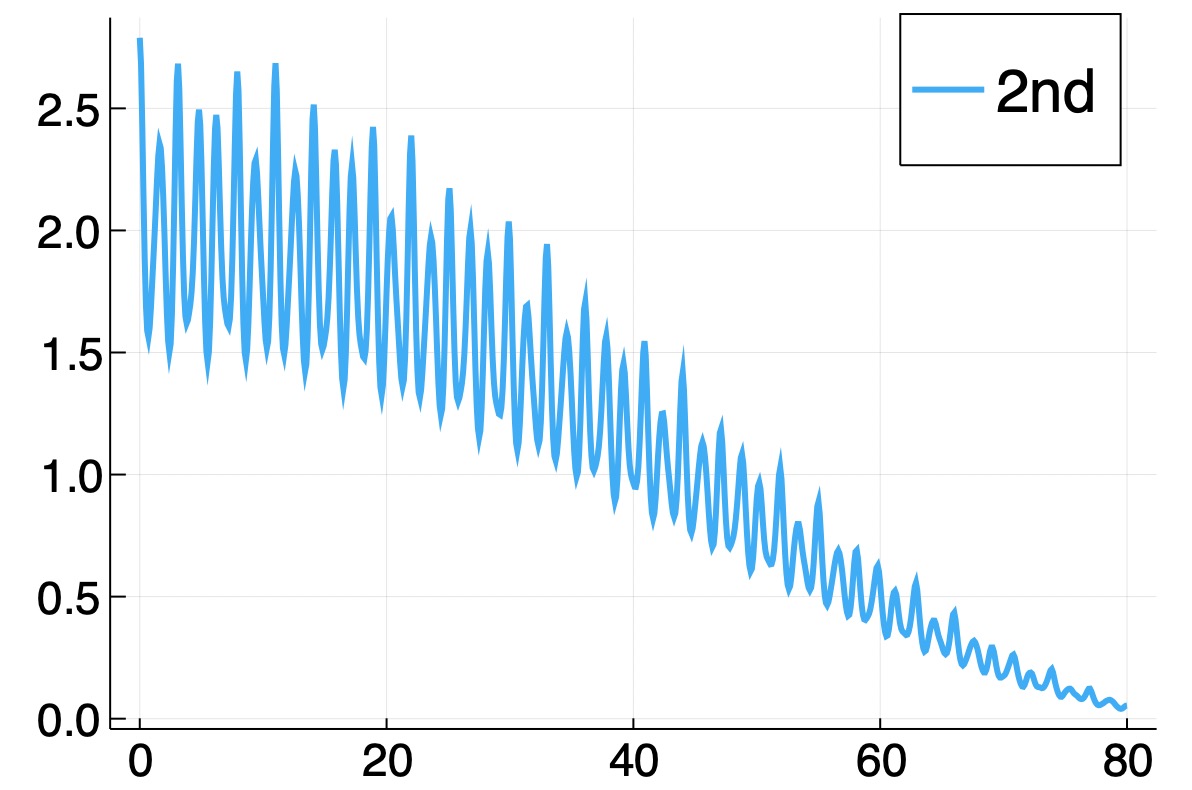}
	\includegraphics[width=4.8cm]{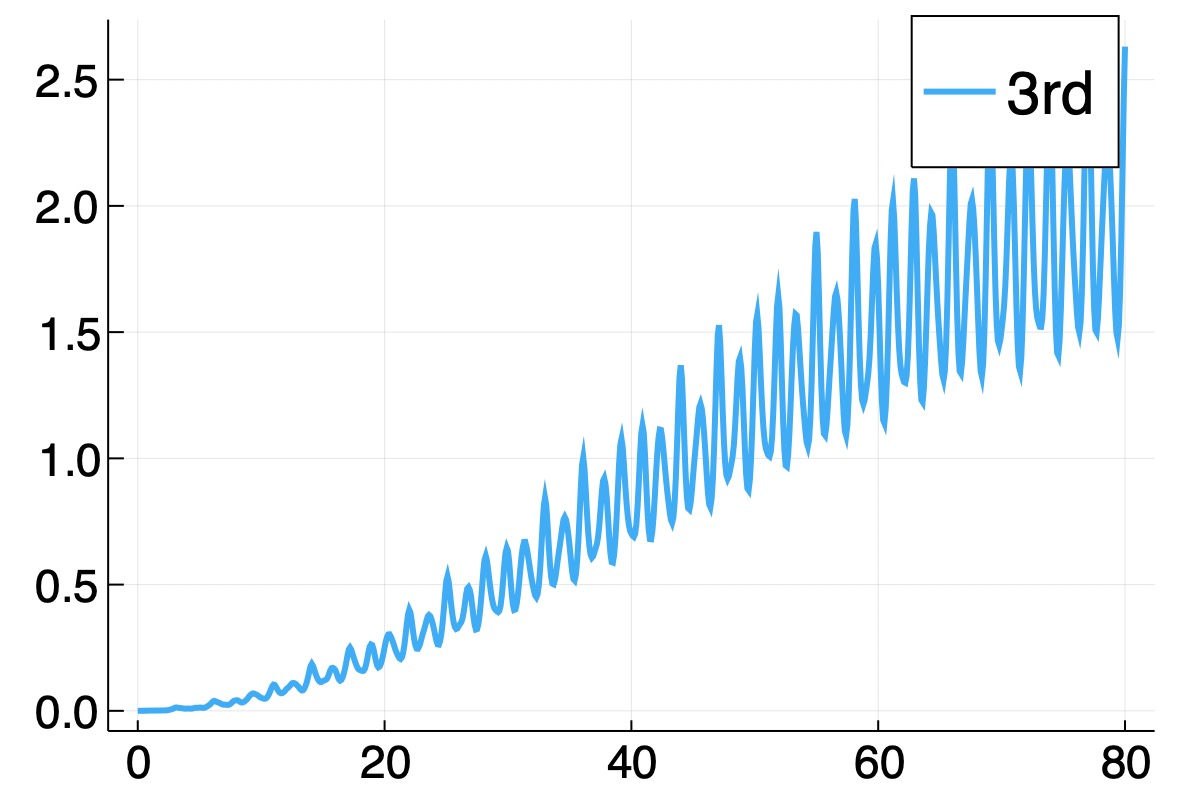}
	\includegraphics[width=4.8cm]{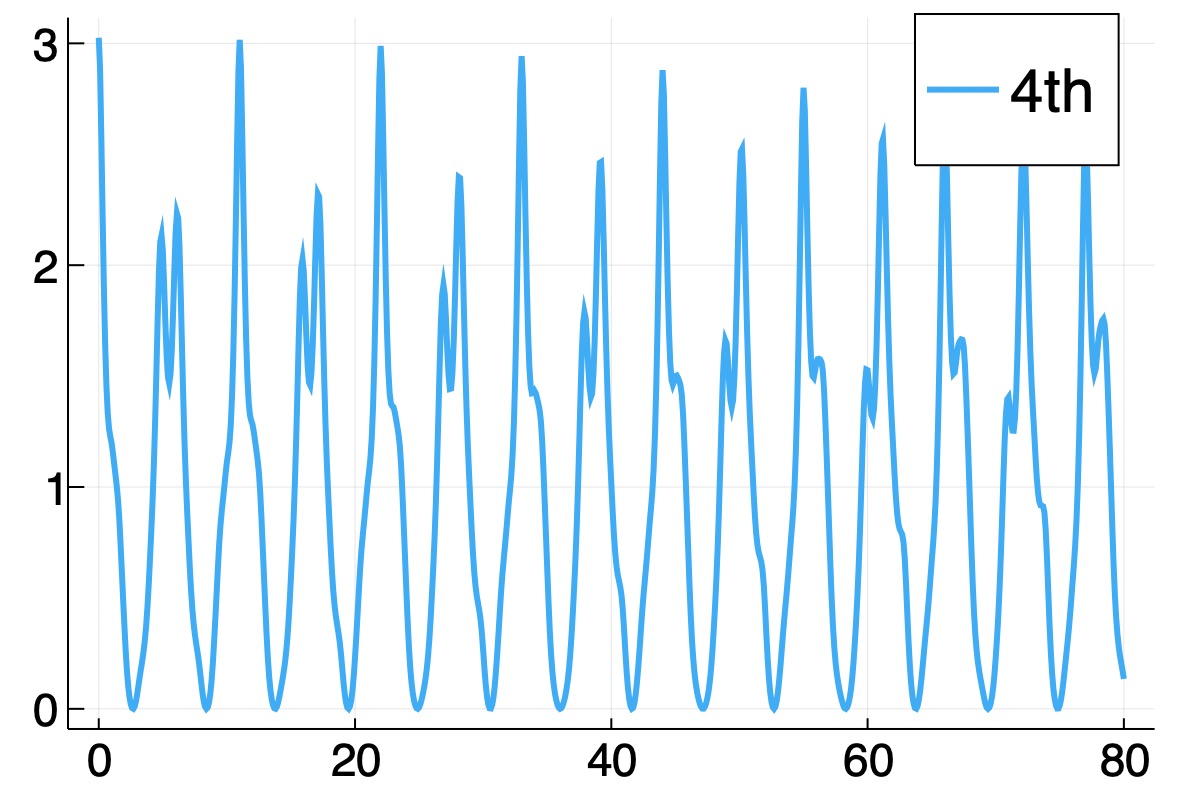}
	\includegraphics[width=4.8cm]{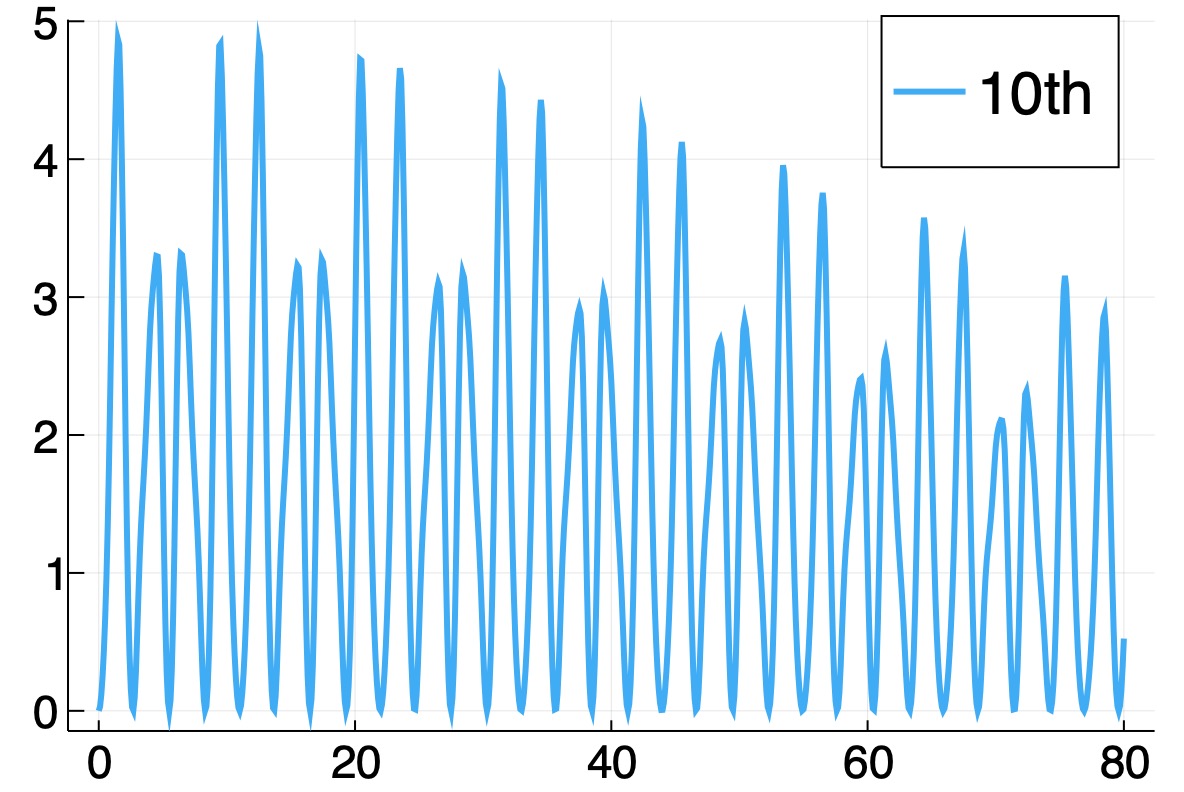}
	\includegraphics[width=4.8cm]{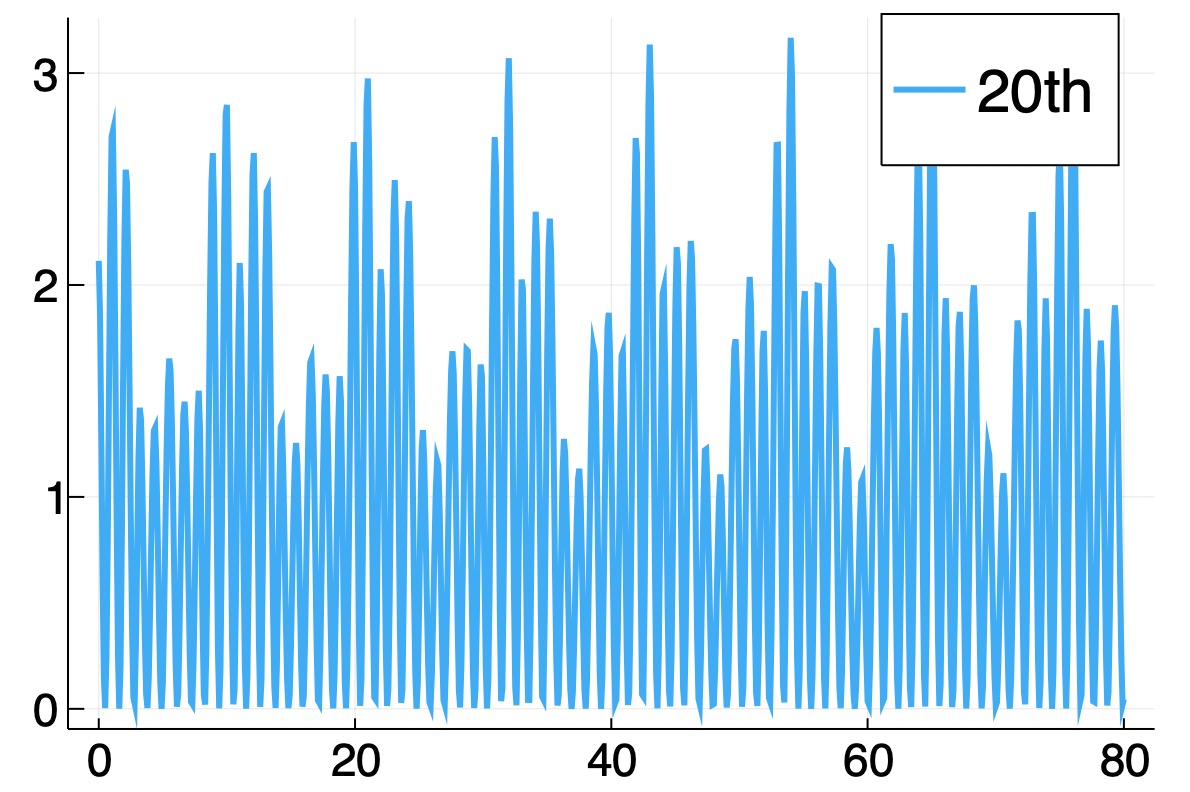}
	\caption{Incommensurate system with $L_1=1$ and $L_2=\pi/2$. The square of norm of the 1st, 2nd, 3rd, 4th, 10th, 20th eigenfunctions with ${\bf k}={\bf 0}$ and $\Ec=2000$.}
	\label{fig:example1:eigenfunctions_a}
\end{figure}

\begin{figure}[htb!]
	\centering
	\includegraphics[width=4.8cm]{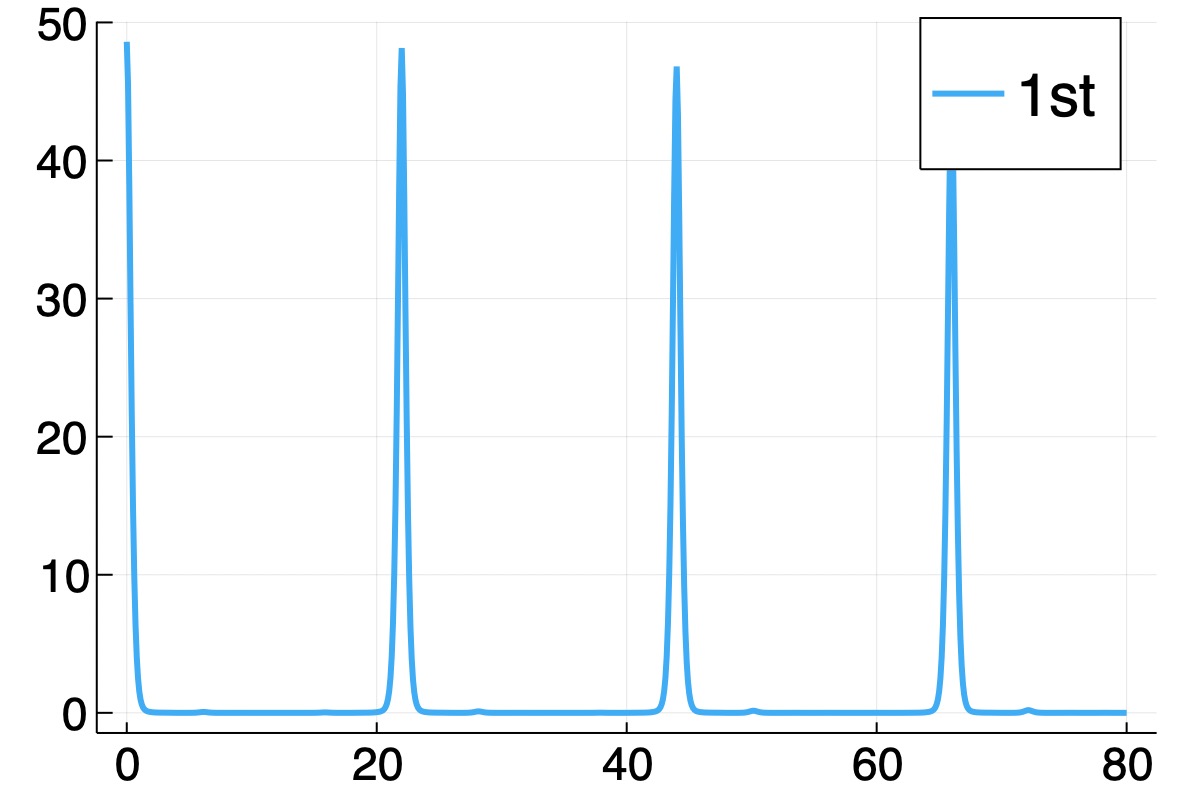}
	\includegraphics[width=4.8cm]{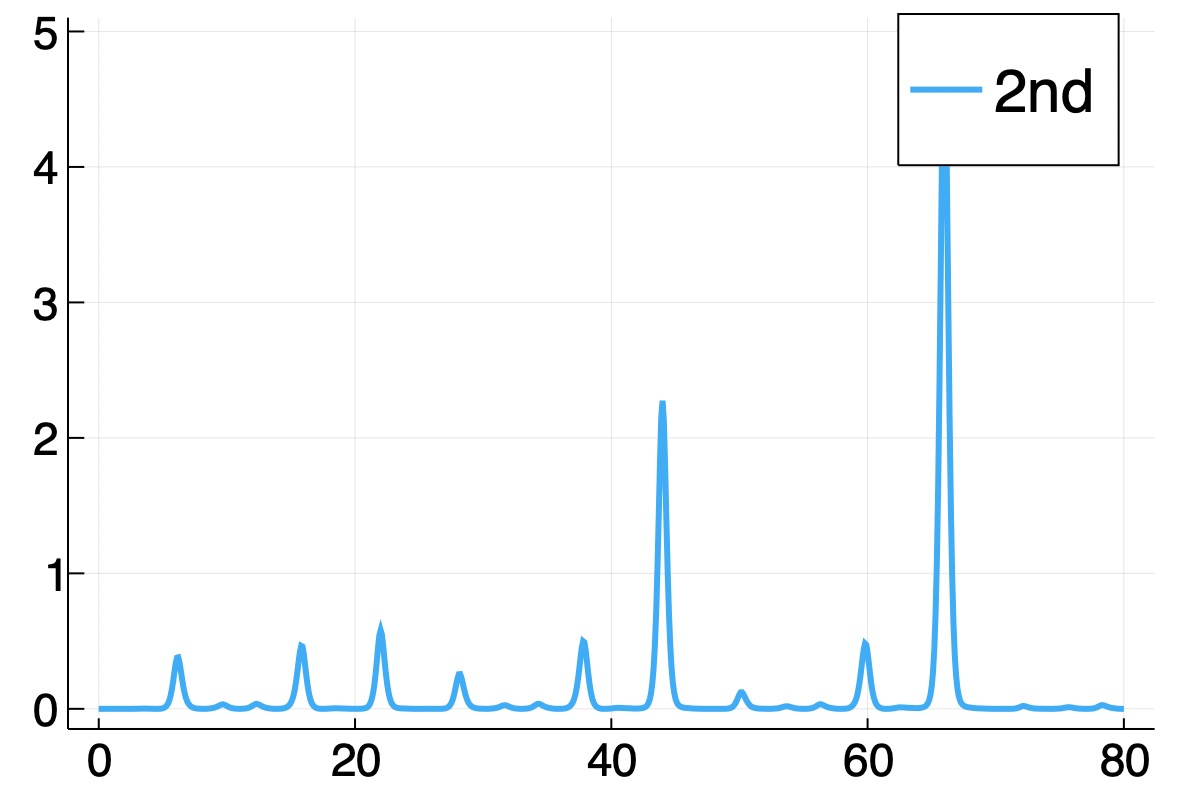}
	\includegraphics[width=4.8cm]{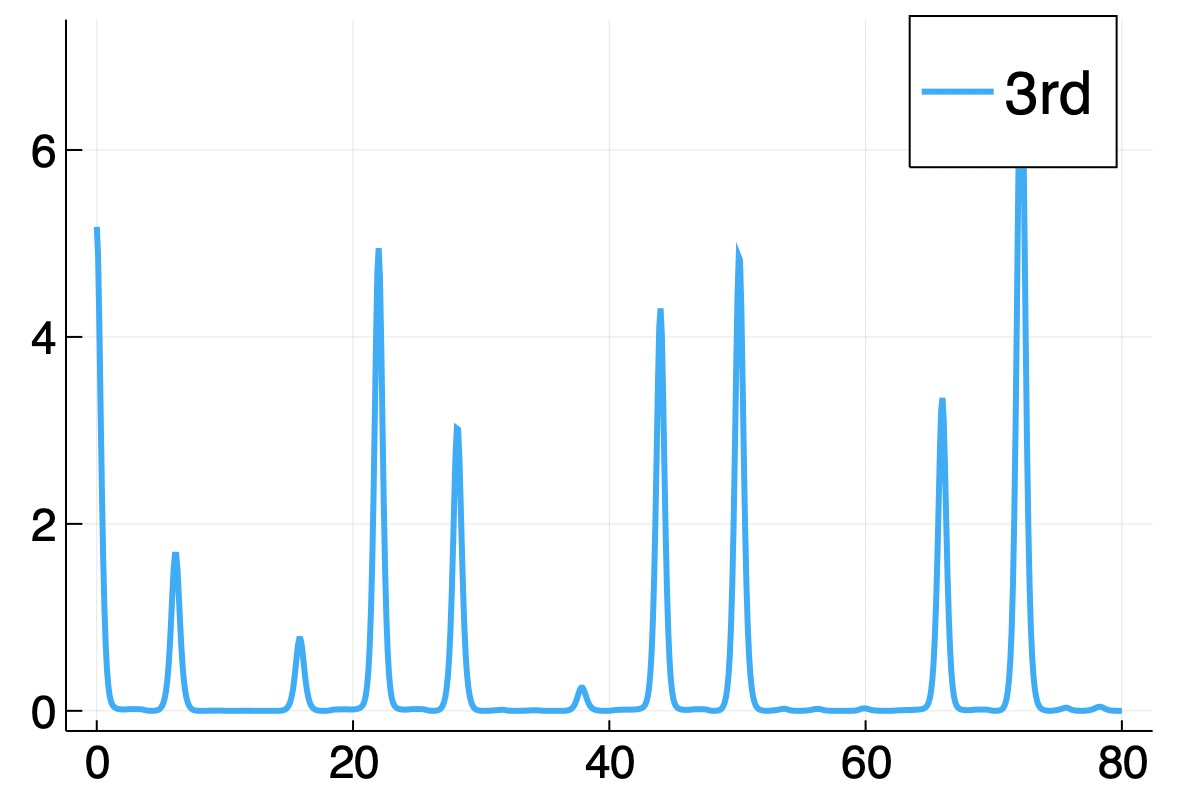}
	\includegraphics[width=4.8cm]{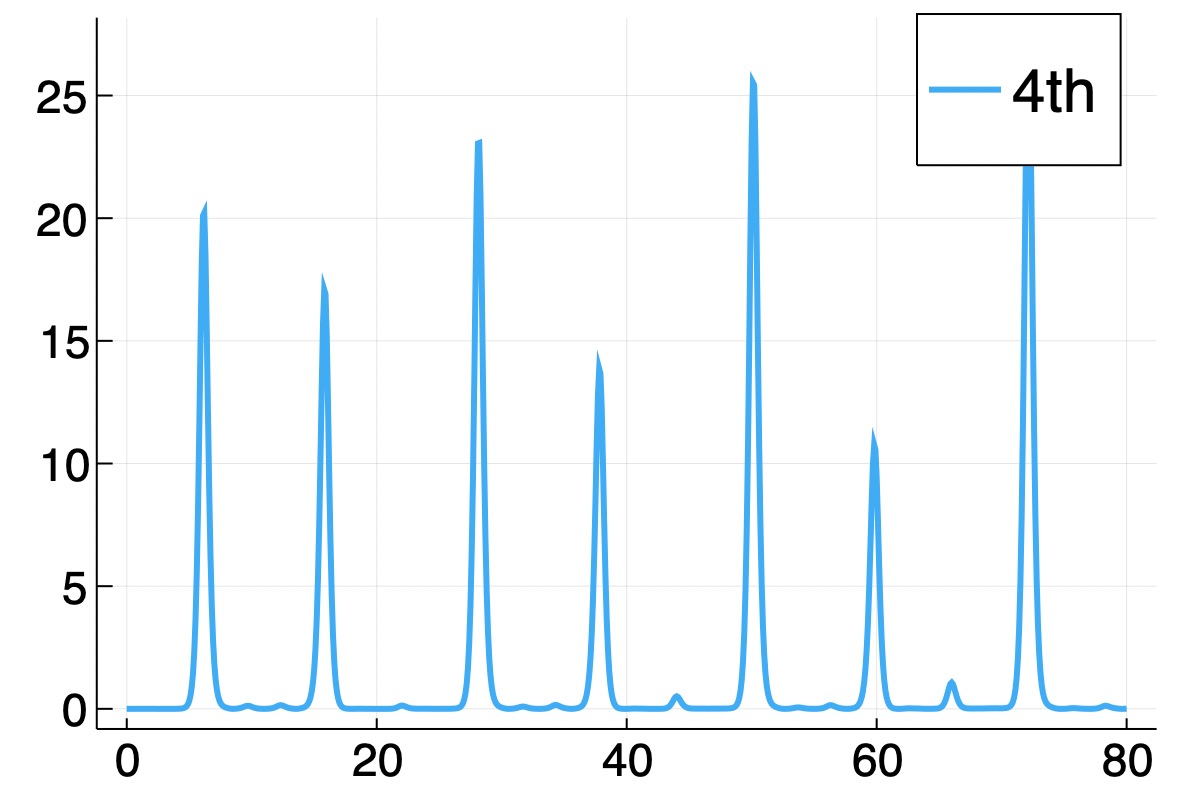}
	\includegraphics[width=4.8cm]{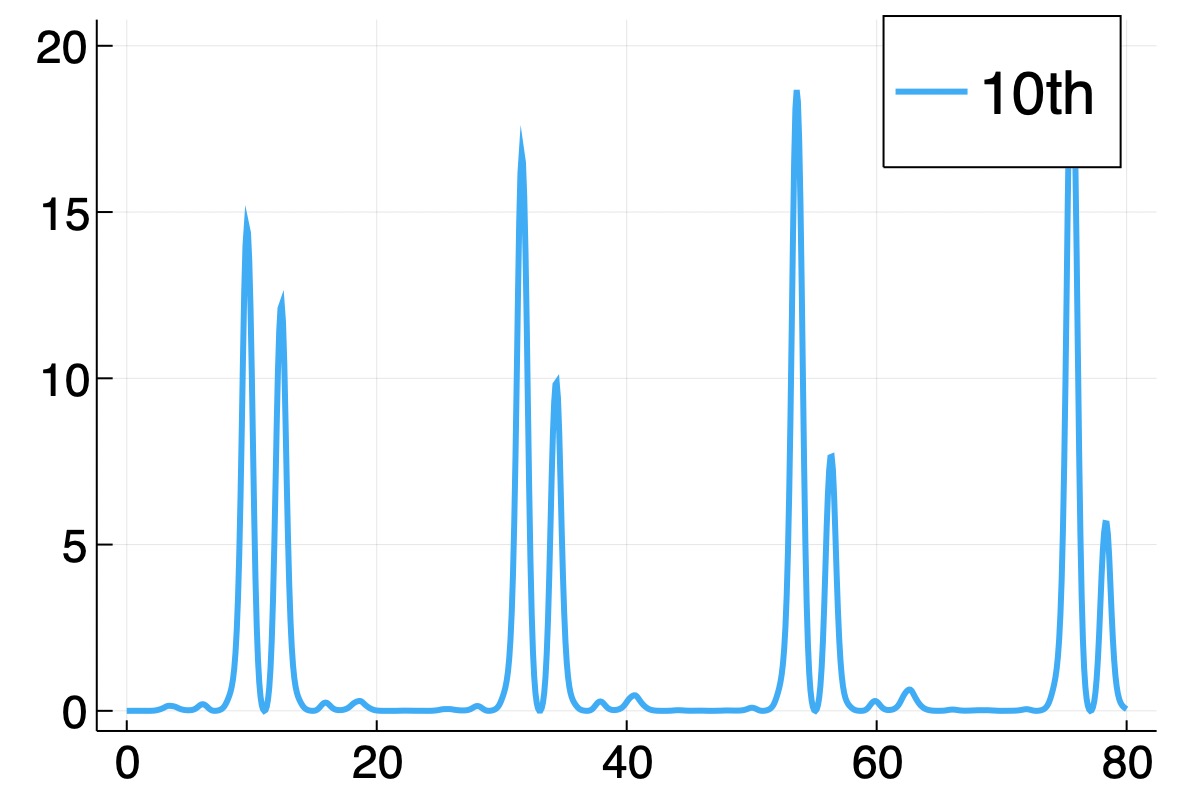}
	\includegraphics[width=4.8cm]{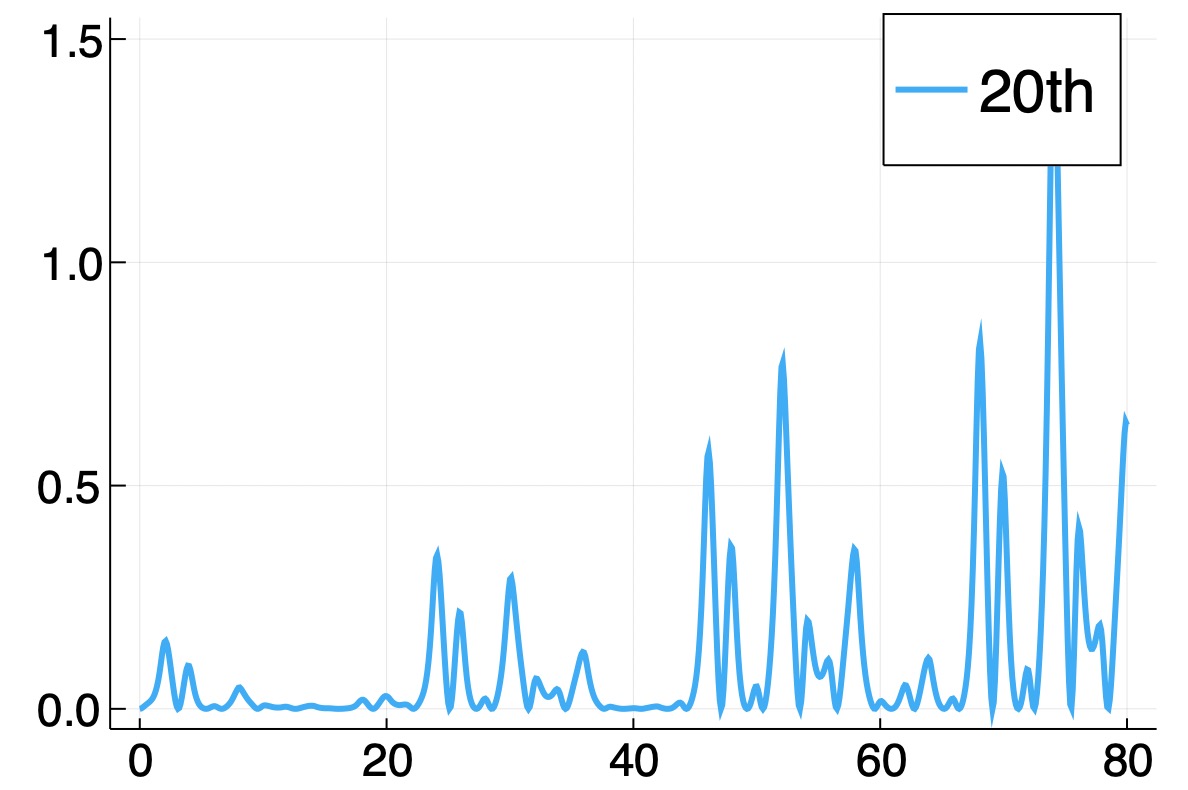}
	\caption{Incommensurate system with $L_1=2$ and $L_2=pi$. The square of norm of the 1st, 2nd, 3rd, 4th, 10th, 20th eigenfunctions with ${\bf k}={\bf 0}$ and $\Ec=2000$.}
	\label{fig:example1:eigenfunctions_b}
\end{figure}

With a  given ${\bf k}={\bf 0}$ and energy cutoff $\Ec=2000$, we plot several eigenfunctions in Fig. \ref{fig:example1:eigenfunctions_a}. 
As a comparison, we also simulate another similar incommensurate system, but with larger lattice constants $L_1=2$ and $L_2=\pi$ for the two atomic chains.
We use the same plane wave discretization and show the eigenfunctions in Fig. \ref{fig:example1:eigenfunctions_b}. 
We observe that in the second system, the eigenfunctions are much more concentrated in some local regions in the second system.

It is worth noting that even in our simple approximation of the atomic potential, the results share similar nature of localization-to-delocalization transitions with many other studies on the 1D bichromatic incommensurate potentials \cite{Li2017i,lschenPRL,Settino2017,Sun2015}, whose major focus is to study the quantum localization. The classical model of the 1D bichromatic incommensurate potentials is the Aubry-Andr\'e (AA) model \cite{AAmodel}, based on the nearest neighbor tight-binding approaches. 
In the AA model, it is predicted that all the electronic states are either localized in the real space or in the reciprocal space (delocalized in the real space), which is determined by the competition of the two potential strength. In other words, there is no mobility edge, i.e., a critical energy separating localized and delocalized energy eigenstates. 
However, according to recent theoretical studies \cite{Li2017i,Settino2017,Sun2015}, and more recent verifications from experiment \cite{lschenPRL}, the mobility edge does exist in the 1D bichromatic incommensurate potentials. 
The reason is that AA model only considers the nearest neighbor interaction, hile recent theoretical models capture more or less continuum nature, by either extending to non-nearest-neighbor hopping within the tight-binding model \cite{Sun2015}, or solving the Hamiltonian semi-continuously through discretizing the real space coordinates plus numerical methods to diagonalize the resulting matrices \cite{Li2017i,Settino2017}.  The numerical schemes applied in these references rely on the commensurate supercell approximations, and cannot be easily extended to more general electronic structure calculations.
Our simulation results in Fig.~\ref{fig:example1:eigenfunctions_b} also manifest such nature. 
For lower eigenstates, the eigenfunction are strongly localize at certain sites, which are separated by few tens of characteristic atomic lengths.
As the eigenvalue goes higher, the corresponding states become more and more delocalized and spread across the lattice.
The continuum nature has been captured by our simulations, which is crucial for the mobility edge to show up in the 1D incommensurate systems.

\begin{figure}[htb!]
	\centering
	\includegraphics[width=14cm]{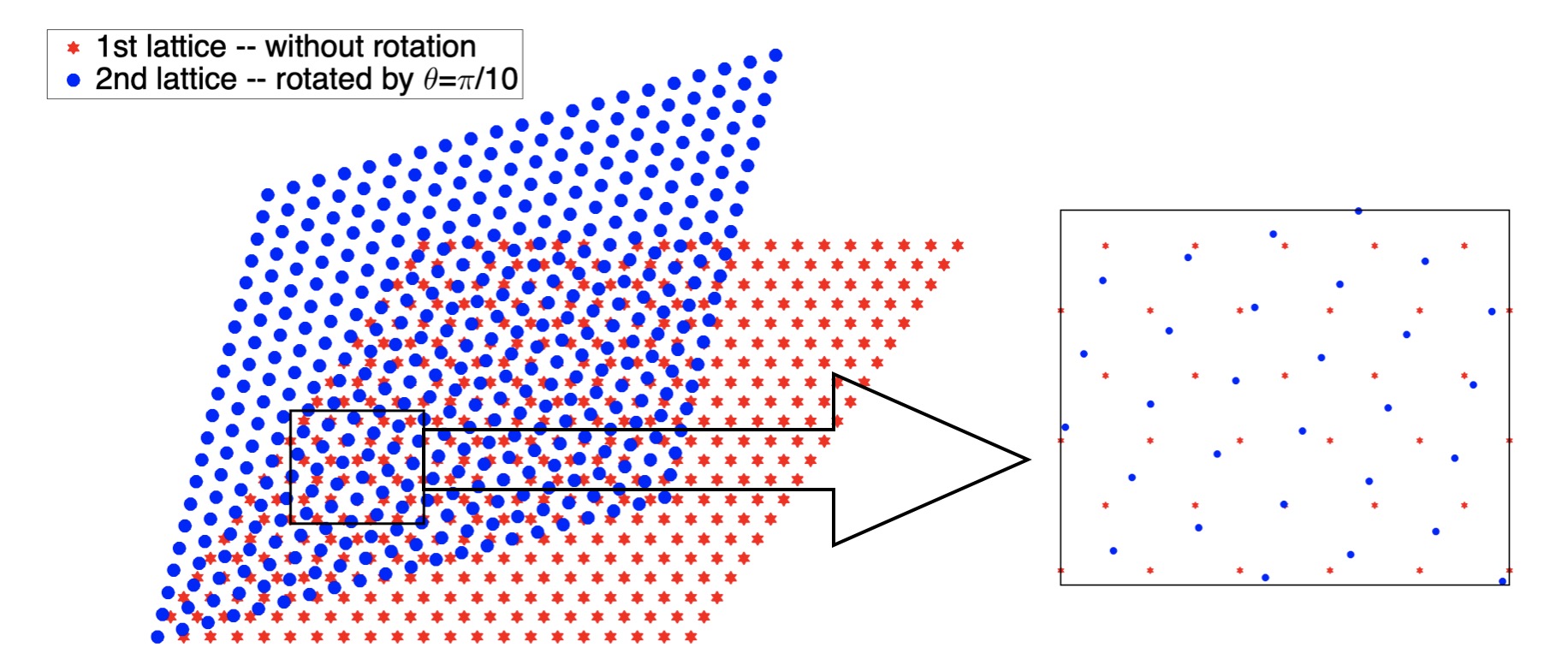}
	\caption{Configuration of the 2D incommensurate layered system, which is from a rotation of two identical triangular lattice.}
	\label{fig:2d_lattice}
\end{figure}

\begin{figure}[!htb]
	\begin{minipage}[t]{0.5\linewidth}
		\centering
		\includegraphics[width=7cm]{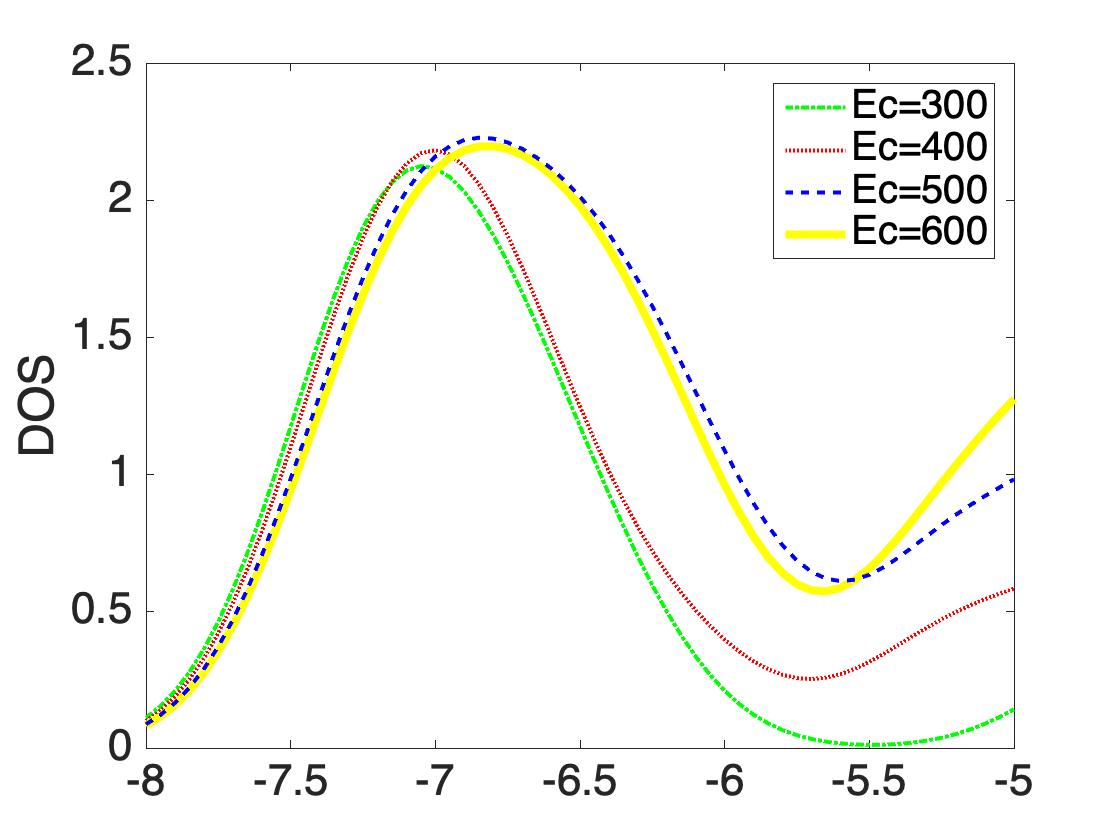}
		\caption{Convergence of $\Dos(\epsilon)$ with respect to energy cutoff $\Ec$ }
		\label{fig:example2:dos:convergence_Ec}
	\end{minipage}
	\hspace{0.3cm}
	\begin{minipage}[t]{0.5\linewidth}
		\centering
		\includegraphics[width=7cm]{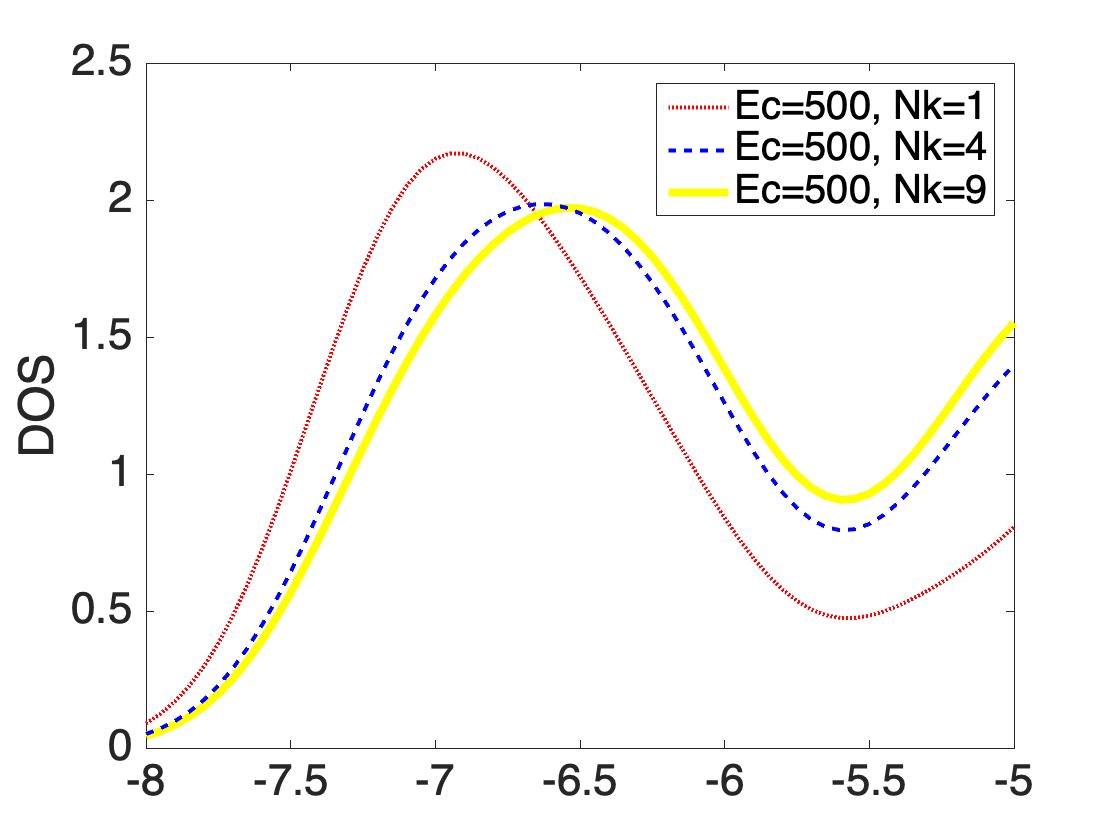}
		\caption{Convergence of $\Dos(\epsilon)$ with respect to ${\bf k}$-point sampling. }
		\label{fig:example2:dos:convergence_k}
	\end{minipage}
\end{figure}

\noindent
{\bf Example 2.}
(two-dimensional sheets with incommensurate rotations).
Consider an incommensurate triangular bilayer, in which one sheet is rotated by $\theta = \pi/10$
with respect to the other  (see Fig. \ref{fig:2d_lattice}).
More precisely, we take $\RL_1=A_1\Z^d$ and  $\RL_2=A_2\Z^d$  with
\begin{eqnarray*}
	A_1=L \cdot \left[ \begin{array}{cc} 1 & \frac{1}{2} \\  0 & \frac{\sqrt{3}}{2} \end{array} \right]
	\qquad{\rm and}\qquad
	A_2=L \cdot \left[ \begin{array}{cc} \cos(\theta) & \cos(\theta+\frac{\pi}{3}) \\  \sin(\theta) & \sin(\theta+\frac{\pi}{3}) \end{array} \right]
\end{eqnarray*}
and lattice constant $L=2.0$.
We solve the eigenvalue problem
\begin{eqnarray}\label{example2}
-\Delta u({\bf r}) + \big( V_1({\bf r}) + V_2({\bf r}) \big)u({\bf r}) = \lambda u({\bf r}) \qquad {\bf r}\in\R^2,
\end{eqnarray}
where $V_1$ and $V_2$ are screened Coulomb potentials with respect to $\RL_1$ and $\RL_2$ respectively.
The convergence of DoS with respect to energy cutoff $\Ec$ and ${\bf k}$-point sampling are shown in
Fig. \ref{fig:example2:dos:convergence_Ec} and \ref{fig:example2:dos:convergence_k}, respectively.
We observe that the convergence with a single ${\bf k}$-point is slow (especially in the high energy window) since we are not able to apply large enough $\Ec$ for the two dimensional systems.
We will investigate more advanced numerical methods in our future works, for the high dimensional problems.
Moreover, we observe that the convergence with multiple ${\bf k}$-points converge much faster.

\begin{figure}[htb!]
	\centering
	\includegraphics[width=5.5cm]{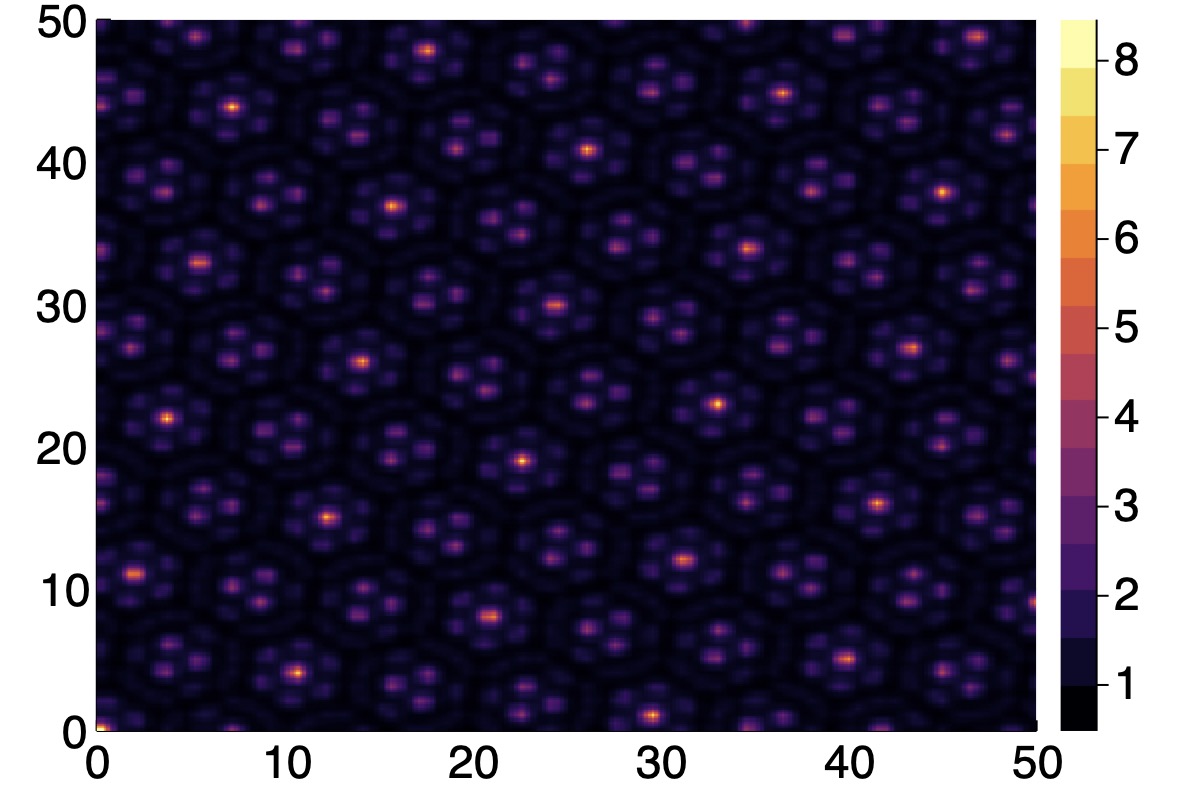}
	\hskip 0.5cm
	\includegraphics[width=5.5cm]{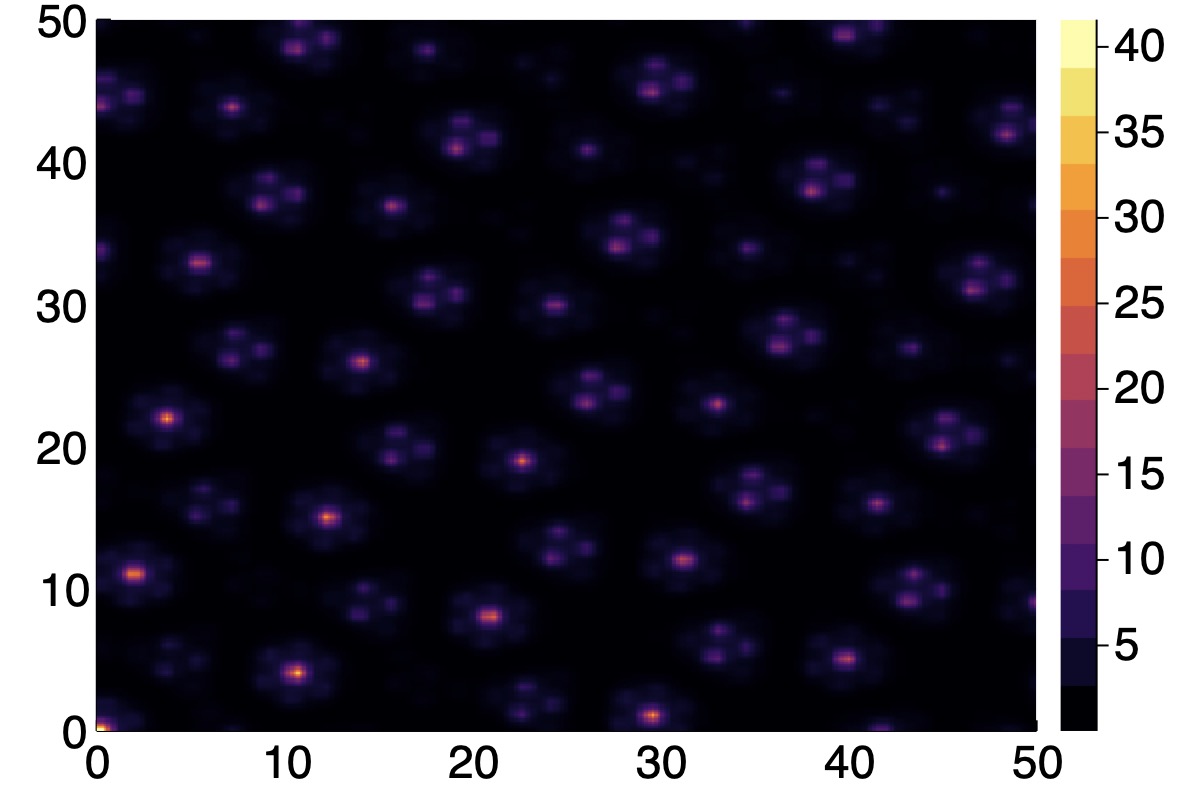}
	\includegraphics[width=5.5cm]{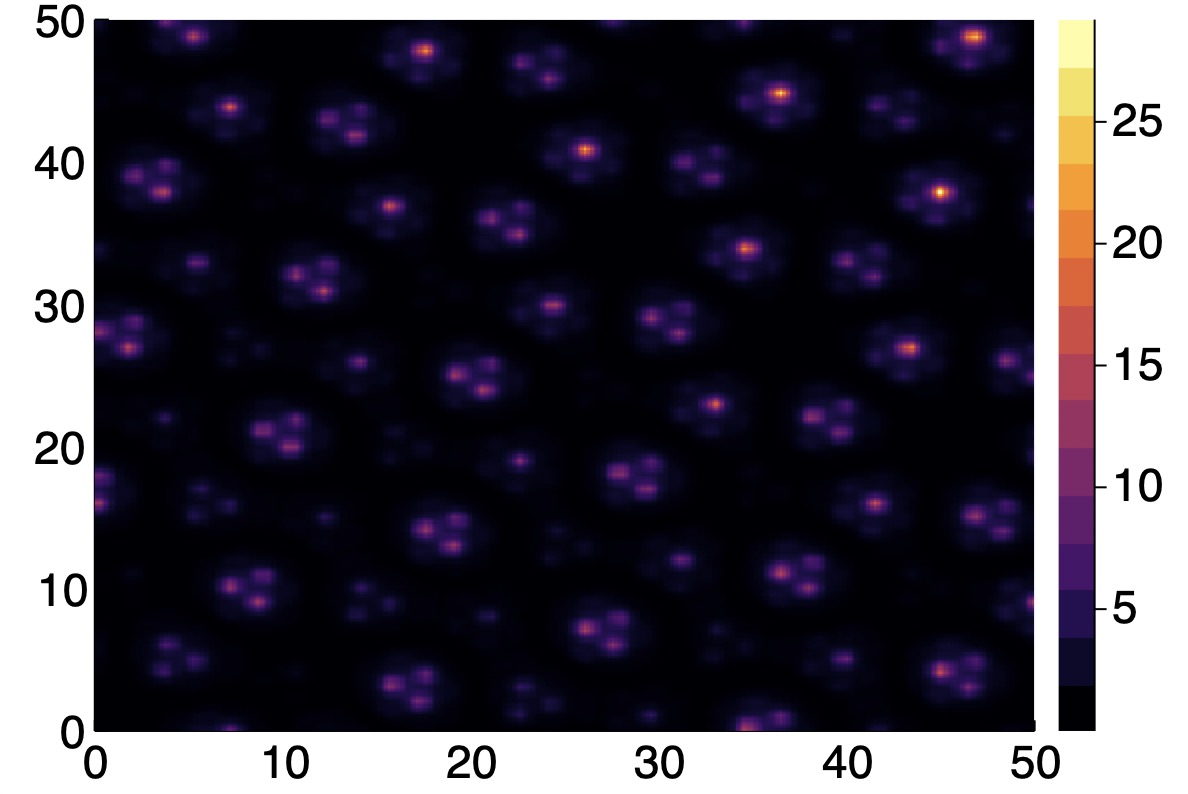}
	\hskip 0.5cm
	\includegraphics[width=5.5cm]{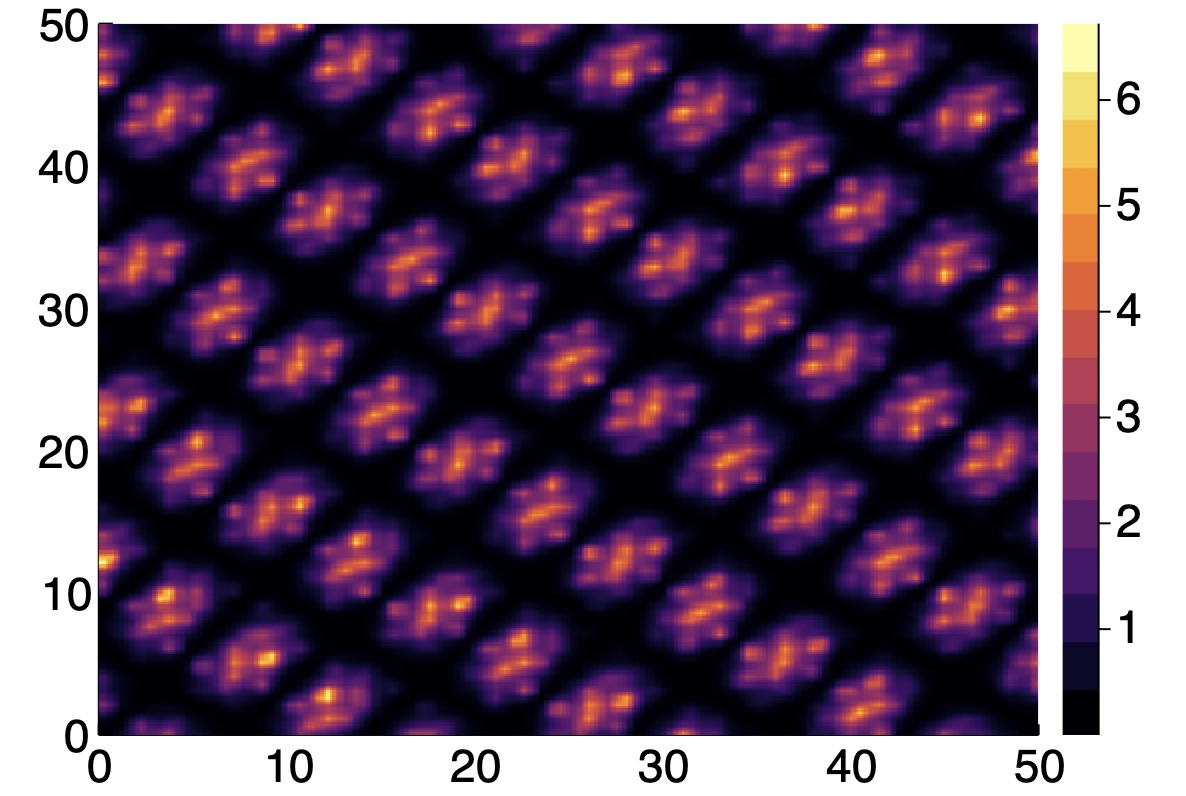}
	\caption{2D incommensurate system with the rotation angle $\theta=\pi/10$. The square of norm of the 1st, 2nd, 3rd, 10th eigenfunctions with ${\bf k}={\bf 0}$ and $\Ec=1000$.}
	\label{fig:example2:eigenfunctions_a}
\end{figure}

\begin{figure}[htb!]
	\centering
	\includegraphics[width=5.5cm]{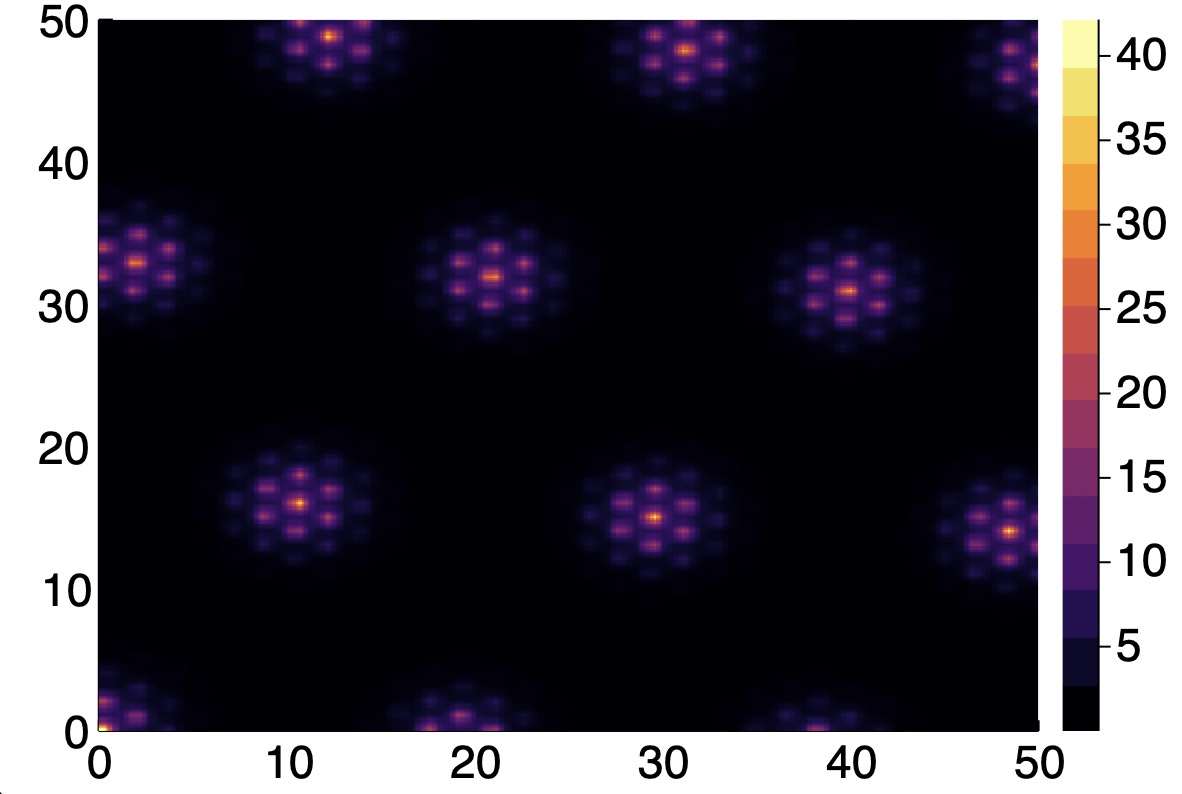}
	\hskip 0.5cm
	\includegraphics[width=5.5cm]{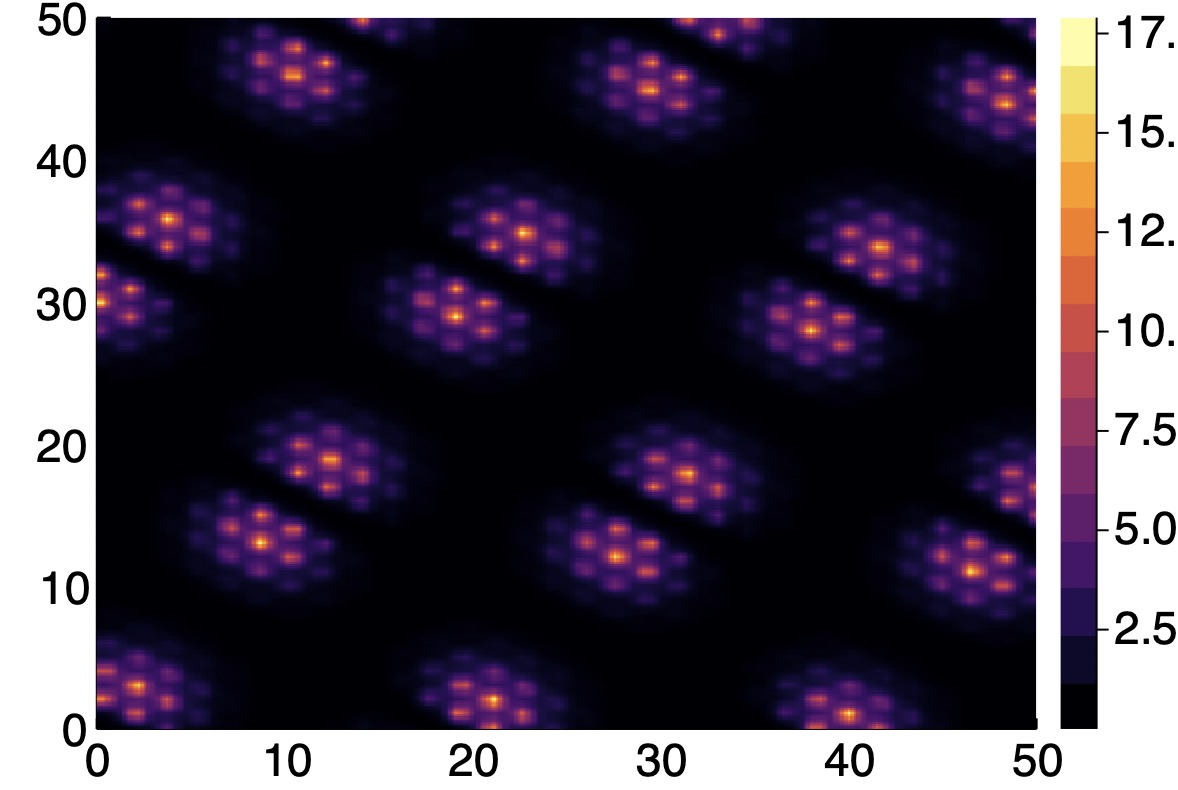}
	\includegraphics[width=5.5cm]{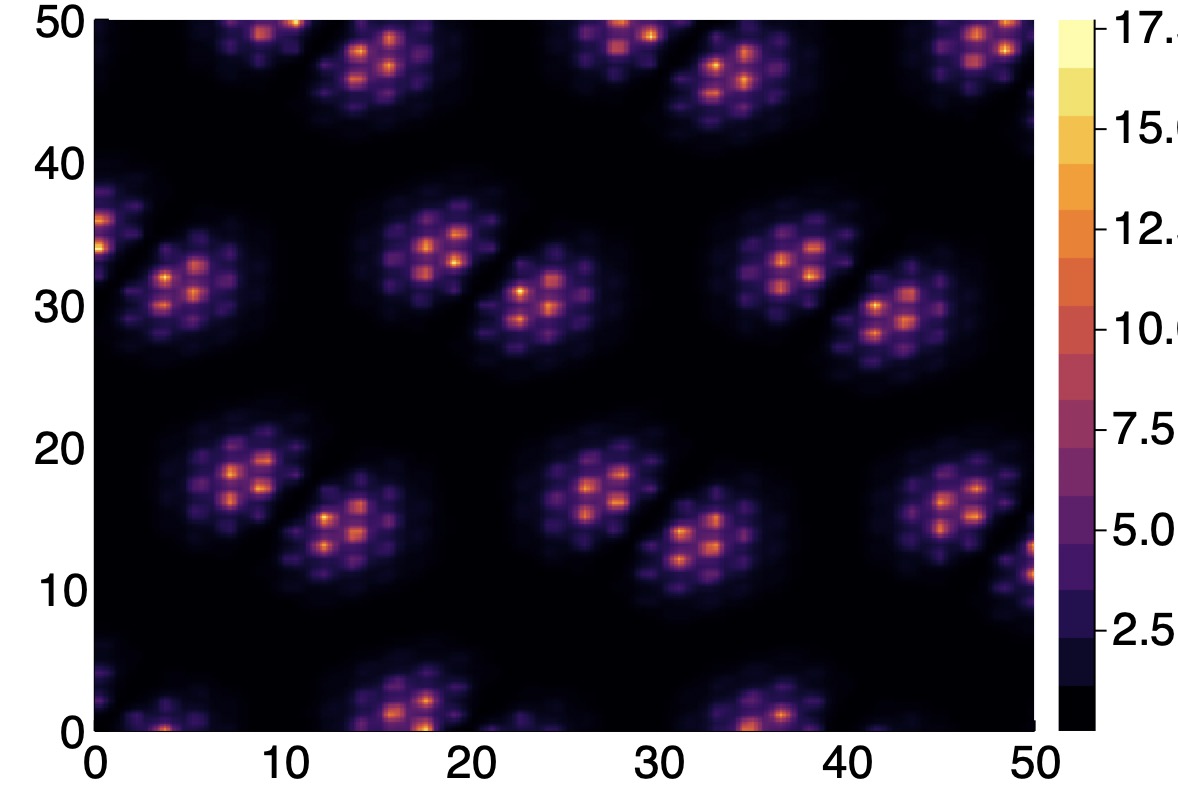}
	\hskip 0.5cm
	\includegraphics[width=5.5cm]{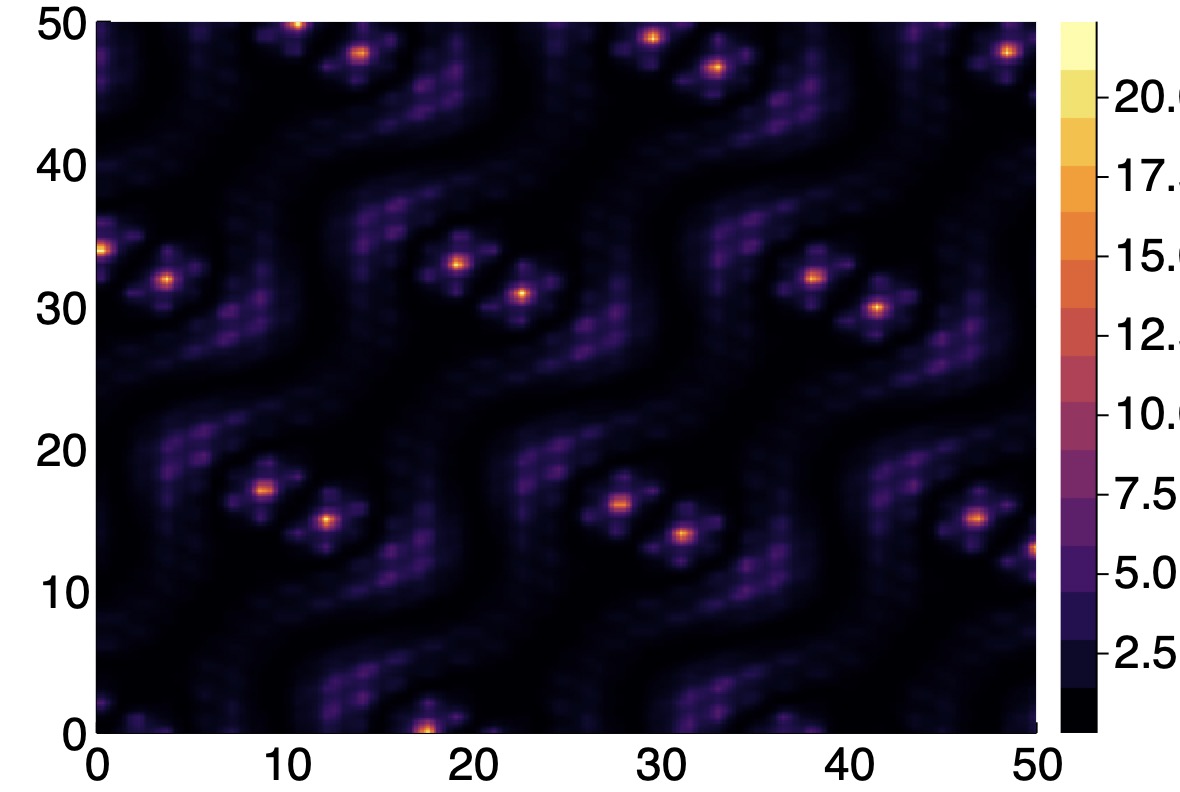}
	\caption{2D incommensurate system with the rotation angle $\theta=\pi/30$. The square of norm of the 1st, 2nd, 3rd, 10th eigenfunctions with ${\bf k}={\bf 0}$ and $\Ec=1000$.}
	\label{fig:example2:eigenfunctions_b}
\end{figure}

With a Gamma point (${\bf k}={\bf 0}$) and energy cutoff $\Ec=1000$, we plot some eigenfunctions of this incommensurate problem in Fig. \ref{fig:example2:eigenfunctions_a}.
For comparison, we also simulate a slightly different incommensurate system with the same lattice constant $L=2$, but a rotation angle $\theta=\pi/30$.
We use the same plane wave discretization and show the eigenfunctions in Fig. \ref{fig:example2:eigenfunctions_b}. 
We observe that the eigenfunctions of the second system are significantly localized. 
Similar to the discussions for the 1D example, we have shown that there are also strong localization effects in 2D incommensurate systems.

\section{Extension to full DFT calculations}
\label{sec:further}
\setcounter{equation}{0}

Without too much difficulties, we can extend our method to full Kohn-Sham DFT calculations,
which is more of practical interest. In this section, we will show that, within the Kohn-Sham DFT framework,
the electron density and total energy (and the variational form) of incommensurate systems
can be naturally expressed by the current plane wave methods.
Furthermore, we will discuss the formulas of DoS and choice of Fermi level.
Our focus of this paper is to give a general framework of solving the quantum problems for incommensurate systems,
so the simulations of real systems will be reserved for our future works.

We consider two periodic lattices $\RL_1,\RL_2\in\R^d$, such that the two lattices $\RL_1$ and $\RL_2$ are incommensurate.
Using the notations in Section \ref{sec:incommensurate}, we denot the unit cells in real space by $\Gamma_1$ and $\Gamma_2$,
the reciprocal lattices by $\RL_1^*$ and $\RL_2^*$, the reciprocal unit cells by $\Gamma_1^*$ and $\Gamma_2^*$.
Moreover, we denote by $Z_1, Z_2\in\Z$ the nuclear charge (or the number of electrons) per unit cell on each lattice.

\subsection{Electron density}
\label{sec:rho}

In the first place, the expression of the electron density must be obtained. Using the plane wave discretizations with a single ${\bf k}$-point (e.g. ${\bf k}=0$)  and an energy cutoff $\Ec$,
the general form of electron density treated in independent-particle theory  can be written as
\begin{eqnarray}\label{density}
\rho^{\bf k}({\bf r}) = \sum_{1\leq j\leq \Nc} f(\lambda_j^{\bf k}) |u_{j,{\bf k}}({\bf r})|^2,
\end{eqnarray}
where $\lambda_j^{\bf k}$ and $u_{j,{\bf k}}$ are the eigenvalue and Kohn-Sham orbital of state $j$, and $f(\lambda_j^{\bf k})$ is the probability of finding an electron in state $j$.
At zero temperature, $f(x) = 2\chi_{(-\infty,\Ef)}(x)$ with $\chi$ the characteristic function, $\Ef$ the Fermi energy and 2 the factor account for spin.
At finite temperature $T$, $f(x) = 2\big(1+\exp((x-\Ef)/(k_{\rm B}T)\big)^{-1}$ with $k_{\rm B}$ the Boltzmann constant.
The choice of Fermi energy will be discussed later in this section.
We mention that the definition of \eqref{density} can be easily generated to multiple ${\bf k}$-points.
For simplicity of presentations, we will omit the subscript ${\bf k}$ in \eqref{density} whenever it is clear from the context.

Using the expression \eqref{eigenfunction}, we have
\begin{eqnarray*}
	\rho({\bf r}) = \sum_{1\leq j\leq \Nc} \sum_{m,n,m',n'}^{\Nc}
	f(\lambda_j) \hat{u}^*_{j}(G_{1m}+G_{2n}) \hat{u}_{j}(G_{1m}+G_{2n}) e^{i(G_{1m}-G_{1m'}+G_{2n}-G_{2n'})\cdot{\bf r}}  .
\end{eqnarray*}
Even though the expression is straightforward, it is not an efficient way to calculate the density,
since finding all the Fourier components of $\rho$ involves a double sum.
To efficiently calculate the electron density, we utilize the fast Fourier transform (FFT) in the higher dimensional formulations.
As discussed in Section \ref{sec:highDim}, we can write $u_j({\bf r})=\tilde{u}_j({\bf r},{\bf r})$
with $\tilde{u}_j({\bf r},{\bf r}')$ the eigenfunction of the higher dimensional problem, which retains the periodicity and has the form of \eqref{eq:eigen-H}.
Therefore, we can construct $\tilde{\rho}({\bf r},{\bf r}'):\R^d\times\R^d\rightarrow\R$
from $\tilde{u}_j({\bf r},{\bf r}')$ (with similar formula as \eqref{density}),
and calculate the electron density by
\begin{eqnarray}\label{rho-dim}
\rho({\bf r}) = \tilde{\rho}({\bf r},{\bf r}) .
\end{eqnarray}
Due to the periodicity, $\tilde{\rho}({\bf r},{\bf r}')$ can be evaluated
in real space (on the grids) and reciprocal spaces (with inverse FFT).
More precisely, we have
\begin{eqnarray}\label{FFT_rho}
\hat{\rho}(G_{1m}+G_{2n}) = \widehat{\tilde{\rho}}(G_{1m},G_{2n}).
\end{eqnarray}
Each Fourier component is unique in the incommensurate systems since $G_{1m}+G_{2n}=G_{1m'}+G_{2n'}$ if and only if $m=m'$ and $n=n'$.

The major advantage of this (higher dimensional) FFT calculations is that $\tilde{\rho}({\bf r},{\bf r}')$
can be directly used to obtain the exchange-correlation term
$\epsilon_{\rm xc}\big( \tilde{\rho}({\bf r},{\bf r}') \big)$ and
$v_{\rm xc}\big( \tilde{\rho}({\bf r},{\bf r}') \big)$  in real space, on the grids.
By performing the inverse FFT transform in higher dimension $\R^d\times\R^d$, we can obtain the Fourier components in the original dimension by
\begin{eqnarray}\label{FFT_xc}
\hat{\epsilon}_{\rm xc}[\rho](G_{1m}+G_{2n}) := \widehat{\epsilon_{\rm xc}\big( \tilde{\rho}\big)}(G_{1m},G_{2n})
\quad{\rm and}\quad
\hat{v}_{\rm xc}[\rho](G_{1m}+G_{2n}) := \widehat{v_{\rm xc}\big( \tilde{\rho}\big)}(G_{1m},G_{2n}).
\quad
\end{eqnarray}
Note that the relation between $\epsilon_{\rm xc}$ and $v_{\rm xc}$ is $v_{\rm xc}(x) = \big(x\epsilon_{\rm xc}(x)\big)'$.

Similarly, we use the Fourier components of electron density to calculate the Hartree energy $E_{\rm H}$
and Fourier components of Hartree potential $\hat{v}_{\rm H}$
\begin{eqnarray}\label{FFT_Hartree}
E_{\rm H} =  \frac{1}{2}\sum_{m,n\neq0}^{\Nc} \frac{\hat{\rho}(G_{1m} + G_{2n})^2}{|G_{1m} + G_{2n}|^2}
\quad{\rm and}\quad
\hat{v}_{\rm H}[\rho](G_{1m}+G_{2n}) := \frac{\hat{\rho}(G_{1m} + G_{2n})^2}{|G_{1m} + G_{2n}|^2}
\quad
\end{eqnarray}
for $m,n\neq 0$.

\subsection{Total energy and Kohn-Sham equations}

With the electron density, Hartree energy, and exchange-correlation energy written in Fourier components,
we can now derive the Kohn-Sham total energy expression for the incommensurate system. 
Denote ${\bf G}_{mn}  =  G_{1m} + G_{2n}$ for $G_{1m}\in\RL^*_1$ and $G_{2n}\in\RL^*_2$.
For a given ${\bf k}$-point and an energy cutoff $\Ec$, the energy  {\it averaged on per unit volume} can be written as
\begin{multline} \label{eq:totalKS}
E^{\bf k}_{\rm  tot}
:= E_{\rm II}  +  \sum_{j} f_j \Bigg\{\sum_{m,n,m'n'}^{\Nc} \hat{u}_{j,{\bf k}}^{*}({\bf G}_{mn})
~ \bigg[ \frac{1}{2} |{\bf k}+{\bf G}_{mn}|^2 \delta_{mm'}\delta_{nn'}
+ V_{1,{\rm ext}}({\bf G}_{mn},{\bf G}_{m'n'})
\\
+ V_{2,{\rm ext}}({\bf G}_{mn},{\bf G}_{m'n'})  \bigg]
~ \hat{u}_{j,{\bf k}}({\bf G}_{m'n'}) \Bigg\}
+ \sum_{mn}^{\Nc} \hat{\epsilon}_{xc}({\bf G}_{mn}) \hat{\rho}({\bf G}_{mn})
+ \frac{1}{2}\sum_{m,n\neq0}^{\Nc} \frac{\hat{\rho}({\bf G}_{mn})^2}{|{\bf G}_{mn}|^2}
\end{multline}
with $U_j^{\bf k}:=\big\{\hat{u}_{j,{\bf k}}({\bf G}_{mn})\big\}_{|G_{1m}|^2+|G_{2n}|^2\leq 2\Ec} \in\R^{\Nc}$ the Fourier components of the $j$th orbital,
$E_{\rm II}$ the nuclei-nuclei interaction energy averaged on per unit volume (we refer to \ref{sec:ewald} for its calculations by using Ewald sum), and
$V_{1,{\rm ext}}$ and $V_{2,{\rm ext}}$ the external potentials generated from two periodic layers $\RL_1$ and $\RL_2$ seperately.
The above expression can be easily generated to multiple ${\bf k}$-points.
For simplicity, we will omit the subscript ${\bf k}$ in \eqref{eq:totalKS} whenever it is clear from the context.

Since $E_{\rm II}$ is a constant with a fixed atomic configuration, the ground state solution of the incommensurate system can be obtained by minimizing the total energy \eqref{eq:totalKS}
with respect to $\big\{U_j^{\bf k}\big\}$ under the orthonormal constraints $U_i^{{\bf k}*}U_j^{\bf k} = \delta_{ij}$.
The variational form of this minimization problem is the {\it discrete} Kohn-Sham equation:
\begin{eqnarray}\label{eq:KS-discrete}
H^{\bf k}[\rho] U_j^{\bf k}  = \lambda_j^{\bf k} U_j^{\bf k}
\end{eqnarray}
with the matrix elements
\begin{multline}\label{hamiltonian_matrix_element}
\qquad
H^{\bf k}[\rho]_{mn,m'n'} = \frac{1}{2} \big|{\bf k}+{\bf G}_{mn}\big|^2 \delta_{mm'}\delta_{nn'}
+ V_{1,{\rm ext}}({\bf G}_{mn},{\bf G}_{m'n'})
+ V_{2,{\rm ext}}({\bf G}_{mn},{\bf G}_{m'n'})
\\
+ \hat{v}_{\rm xc}[\rho]({\bf G}_{mn}-{\bf G}_{m'n'}) + \hat{v}_{\rm H}[\rho]({\bf G}_{mn}-{\bf G}_{m'n'})
\qquad\qquad
\end{multline}
for $|G_{1m}|^2+|G_{2n}|^2\leq 2\Ec$ and $|G_{1m'}|^2+|G_{2n'}|^2\leq 2\Ec$.

Eq. \eqref{eq:KS-discrete} is a nonlinear eigenvalue problem due to the dependence of  $H^{\bf k}[\rho]$ on $\rho$,
and hence the eigenvectors $U_j^{\bf k}$.
We shall resort to the self-consistent field (SCF) iterations to solve this problem.
In each step of the iteration, a linear eigenvalue problem needs to be solved to obtain the trival electron density for next step.

The eigenvalue problem of full DFT calculation has similar form of \eqref{eq:eigen-H-matrix},
with two slight differences.
First, in the pseudopotential approaches \cite{martin04}, the external potential $V_{i,{\rm ext}}$ includes
a local potential $V_{i,{\rm loc}}$ and a nonlocal operator $V_{i,{\rm nl}}$ for $i=1,2$.
Since both $V_{i,{\rm loc}}$ and $V_{i,{\rm nl}}$ are translation invariant with respect to $\RL_i$ for $i=1,2$,
we can compute the matrix elements by
\begin{eqnarray}\label{H_ext}
V_{1,{\rm ext}}({\bf G}_{mn},{\bf G}_{m'n'}) = \delta_{nn'} \Big( \hat{V}_{1,{\rm loc}}(G_{1m}-G_{1m'})
+ \intRd e^{-i{G_{1m}}{\bf r}} V_{1,{\rm nl}} e^{i{G_{1m'}}{\bf r}} \dd{\bf r}  \Big)
\end{eqnarray}
and the same formula for $V_{2,{\rm ext}}({\bf G}_{mn},{\bf G}_{m'n'}) $.
%
Second, the Hartree and exchange-correlation potentials in \eqref{hamiltonian_matrix_element} are neither
periodic with respect to $\RL_1$ nor $\RL_2$.
Despite the two additional terms compared with \eqref{eq:eigen-H-matrix},
we have from the discussions in Section \ref{sec:rho}, \eqref{hamiltonian_matrix_element} and \eqref{H_ext} that
these two terms are periodic with respect to $\widetilde{\RL}=\RL_1\times\RL_2$ in the higher dimension.
Therefore, we can still apply the same framework in Section \ref{sec:planewave} for solving the linear problem
in each SCF iteration step and then the full Kohn-Sham equations.

\subsection{Fermi level determination}
\label{sec:fermi}

Based on the previous analysis, every term in \eqref{eq:totalKS} can be expressed in the current framework,  which enables us to solve the eigenvalues $\{\lambda_j^{\bf k}\}$. 
Having the eigenvalues $\{\lambda_j^{\bf k}\}$ alone is not enough to describe the total energy of the incommensurate system. 
A crucial step to connect the eigenvalues to the total energy is the definition and determination of the Fermi level.

To achieve this, we need to more properly scale the DoS associated with some normalization volume in real space.
For simplicity of presentations, we focus on the simulations with a single ${\bf k}$-point with an energy cutoff $\Ec$.
The extension to multiple ${\bf k}$-points will be straightforwd.
Then the set of plane wave vectors used for one matrix eigenvalue problem (e.g. \eqref{eq:KS-discrete}) is $\big\{{\bf k}+{\bf G}_{mn}\big\}_{|G_{1m}|^2+|G_{2n}|^2\leq 2\Ec}$.

The determination of the Fermi level in the incommensurate system can be better understood when we make some analog to the periodic case. Let us first take the first lattice as a reference for better illustration. Later on when we get to the final definition, it will be independent of the reference lattice. 
With large enough cutoffs, the wavevector set generates uniformly distributed ${\bf k}$-points in reciprocal space. Now we divide the reciprocal space by $\Gamma_1^*$ and take a close look at the first cell $\Gamma_1^{*0}$ that contains the origin. 
Given the form of $\big\{{\bf k}+{\bf G}_{mn}\big\}$, each wavevector in $\Gamma_1^{*0}$ can find a set of replica in all other reciprocal unit cells through shifting by $G_{1m}\in\RL_1^*$. This is very similar to the ${\bf k}$-point sampling of the periodic system. 
The difference is that the ${\bf k}$-points within a reciprocal unit cell are now correlated rather than independent as in the case of periodic systems.
Similar to periodic systems, the scaled DoS with respect to $\RL_1$ is defined as
\begin{eqnarray}\label{eq:dosa1d}
\Dos^{\bf k}_1(\epsilon) = \frac{1}{N_1} \sum_j \delta (\epsilon-\lambda_j^{\bf k}) ,
\end{eqnarray}
where $N_1$ is the number of wavevectors appearing in $\Gamma_1^{*0}$. 
And in real space, such sampling corresponds to a normalization volume of $N_1$ unit cells\footnote{
	In the periodic case, the normalization volume is constructed together with the periodic boundary condition in which the wavefunction can be unambiguously defined within such supercell. 
	In our framework for the incommensurate system, the wavefunction cannot, and more importantly, do not need to explicitly resort to the periodic boundary condition to be clearly defined. 
	However, the correspondence between the ${\bf k}$-point sampling and the normalization volume is very helpful in defining the Fermi level, which we just borrow from the periodic case.
}$^,$\footnote{
	One might feel a little bit uncomfortable that a finite normalization volume could contradict with the incommensurate nature. This can be relieved by thinking of the limit of dense ${\bf k}$-point sampling, when spacing of ${\bf k}$-points goes to zero implying the normalization volume goes to infinite, which is compatible with the incommensurate nature. 
	Essentially, the normalization volume is a parameter underpinned by the ${\bf k}$-point sampling. And it is of more practical importance to examine if the current ${\bf k}$-point sampling gives converged spectrum when performing the DFT total energy calculations for the incommensurate system (and periodic system as well!).
}. 
When $N_1$ is large enough, the average number of electrons in the unit cell of $\Gamma_1$ is well defined, which is $|\Gamma_1|(\frac{Z_1}{|\Gamma_1|}+\frac{Z_2}{|\Gamma_2|})$.
The Fermi level $\Ef$ is then determined by filling these electrons to the scaled DoS:
\begin{eqnarray}\label{eq:fermia1d}
|\Gamma_1|\Big(\frac{Z_1}{|\Gamma_1|}+\frac{Z_2}{|\Gamma_2|}\Big) = \int_{\R} \Dos^{\bf k}_1(\epsilon) f(\Ef,\epsilon) \dd\epsilon ,
\end{eqnarray}
where $f(\Ef, \cdot) = 2\chi_{(-\infty,\Ef)}(\cdot)$ at zero temperature
and $f(\Ef, \cdot) = 2\big(1+\exp((\cdot-\Ef)/(k_{\rm B}T)\big)^{-1}$ at finite temperature $T$.

We can repeat the above process using the 2nd lattice as reference and obtain the scaled DoS
\begin{eqnarray}\label{eq:dosb1d}
\Dos^{\bf k}_2(\epsilon) = \frac{1}{N_2} \sum_j \delta (\epsilon-\lambda_j^{\bf k}) ,
\end{eqnarray}
with $N_2$ is the number of wavevectors appearing in $\Gamma_2^{*0}$.
Then the Fermi level $\Ef$ can be obtained by solving
\begin{eqnarray}\label{eq:fermib1d}
|\Gamma_2|\Big(\frac{Z_1}{|\Gamma_1|}+\frac{Z_2}{|\Gamma_2|}\Big) = \int_{\R} \Dos^{\bf k}_2(\epsilon) f(\Ef,\epsilon) \dd\epsilon .
\end{eqnarray}

We shall see that the definition of the Fermi level is independent of the lattice we choose as the reference, i.e. solving \eqref{eq:fermia1d} and \eqref{eq:fermib1d} are equivalent.
We require that the ${\bf k}$-points are dense and uniformly distributed in reciprocal space, which is guaranteed by the ergodicity and large cutoffs. 
Under this condition, the number of ${\bf k}$-points in a region is proportional to the volume $\frac{N_1}{N_2}=\frac{|\Gamma_1^*|}{|\Gamma_2^*|}=\frac{|\Gamma_2|}{|\Gamma_1|}$.
Therefore, we can unify the definiton of the DoS by a scaling with respect to the unit volume:
\begin{eqnarray}\label{eq:dos1d}
\Dos^{\bf k}(\epsilon) := \frac{1}{|\Gamma_2|} \Dos^{\bf k}_1(\epsilon)  = \frac{1}{|\Gamma_1|} \Dos^{\bf k}_2(\epsilon)
= \frac{1}{\overline{N}} \sum_j \delta (\epsilon-\lambda_j^{\bf k})  ,
\end{eqnarray}
where $\overline{N}:=\frac{N_1}{|\Gamma_2|}=\frac{N_2}{|\Gamma_1|}$.
Note that in the large $\Ec$ limit, $\overline{N}$ is proportional to $\sqrt{\Nc}$ and hence $\sqrt{\Ec}$, which also explains the choice of prefactor used in definition \eqref{def:dos}.
Then the Fermi level $\Ef$ can be defined in a unified way by solving
\begin{eqnarray}\label{eq:fermi1d}
\frac{Z_1}{|\Gamma_1|}+\frac{Z_2}{|\Gamma_2|} = \int_{\R} \Dos^{\bf k}(\epsilon) f(\Ef,\epsilon) \dd\epsilon ,
\end{eqnarray}
where the left hand side is the number of electrons per unit volume.
Now the definitions of the DoS and Fermi level are independent of the lattice we choose to work on.

In the previous subsections, we restrict ourself to the single ${\bf k}$-point calculations, which in principle could achieve convergence with extremely large energy cutoff $\Ec$ due to ergodicity in Lemma \ref{lemma:ergodic}. 
However it is generally computationally inefficient to sample the reciprocal space only with one single ${\bf k}$-point, and we have shown (in the numerical tests in Section \ref{sec:numerics}) that using multiple ${\bf k}$-points could accelerate the convergence significantly. 
The scaling of DoS and the determination of the Fermi level in the case of multiple ${\bf k}$-points are essentially the same as Eqs.~\eqref{eq:dos1d} and \eqref{eq:fermi1d},
only that the DoS is further scaled by a prefactor $\frac{1}{N_k}$ with $N_k$ the number of ${\bf k}$-points in a uniform sampling
\begin{eqnarray*}
\Dos(\epsilon) = \frac{1}{N_k} \sum_{\bf k} \Dos^{\bf k}(\epsilon) .
\end{eqnarray*}
A more detailed study on the ${\bf k}$-point sampling technique and  investigation on the convergence with respect to ${\bf k}$-point sampling will be presented in our future work.

Combining Section 5.1 to Section 5.4, we see that our plane wave methods can be extended to full Kohn-Sham DFT calculations, which paves the way for the future study of more realistic incommensurate systems.

\section{Conclusions}
\label{sec:conclusion}
\setcounter{equation}{0}

In this paper, we propose a plane wave framework for the electronic structure related eigenvalue problems of the incommensurate systems. 
The ergodicity emerging in the incommensurate structures give rise many unique features compared to the periodic systems. 
Our methods can also be extended to full Kohn-Sham DFT calculations of the real systems. 
In principle, the algorithm and theory developed in this paper can be extended to incommensurate systems with $p>2$ layers, 
but with fast growing computational complexity with respect to $p$, which calls for advanced numerical techniques. 
Moreover, the convergence rates of DoS with respect to $\Ec$ and related ${\bf k}$-point sampling strategies will be addressed in our future studies.

\appendix
\renewcommand\thesection{\appendixname~\Alph{section}}

\section{Ewald sum of incommensurate systems}
\label{sec:ewald}

\renewcommand{\theequation}{A.\arabic{equation}}
\renewcommand{\thetheorem}{A.\arabic{theorem}}
\renewcommand{\thelemma}{A.\arabic{lemma}}
\renewcommand{\theproposition}{A.\arabic{proposition}}
\renewcommand{\thealgorithm}{A.\arabic{algorithm}}
\renewcommand{\theremark}{A.\arabic{remark}}
\setcounter{equation}{0}

To compute the force of each atom, we also need to calculate the nuclei-nuclei interaction $E_{\rm II}$ in \eqref{eq:totalKS}.
This can be calculated by Ewald sum. 

We mention that systems described in the framework of this paper do not have well-defined $\gamma_{\rm Ewald}$, 
since the distance between the two layers in the $(d+1)$th direction is neglected.
This gives rise to the problem that two nuclei could get arbitrarily close and cause a blow-up in the Coulomb energy.
The systems of more practical interests are two periodic lattices paralleled to each other and separated by some distance in the $(d+1)$th direction.
Our following discussion relates to this practical scenario, and we assume that the displacement (in the $(d+1)$th direction) between the two layers is ${\bf t}\in\R^d$.

The key idea of the Ewald sum is to separate the lattice Coulomb sum into two parts, one in the real space and one in the reciprocal space, and both parts could converge rather quickly. 
From the discussion below we can see that the incommensurate structures present no difficulty in evaluating the Ewald energy.
The Ewald sum within each periodic lattices can be readily calculated by using standard techniques, 
and the extra contribution to be considered comes from the interlayer sum.

For simplicity of presentations, we take the mono atomic lattices as example to demonstrate the calculation of the interlayer contribution.
Specifically, we calculate the average Coulomb interaction of an atom in second lattice with all atoms in first lattice.
The real space sum (energy per unit volume) can be written as
\begin{eqnarray}
E_{\rm int,r} = \frac{Z_{1}Z_{2}}{|\Gamma_{1}||\Gamma_{2}|}
\int_{\Gamma_{1}} \sum_{{\bf R}_{1n}\in\RL_1}\frac{{\rm erfc}(|{\bf r}+{\bf t}-{\bf R}_{1n}|)}{|{\bf r}+{\bf t}-{\bf R}_{1n}|}{\rm d}{\bf r}  ,
\end{eqnarray}
where the summation goes over the lattice sites of the first lattice. 
As a direct consequence of erogdicity in Lemma \ref{lemma:ergodic},
the projection of all atoms of the second lattice onto the first lattice will result in equal likelihood of any positions within the unit cell of first lattice. 
To average this part of energy, we perform the integral in the domain with the geometry of the unit cell of first lattice.

Similarly, the reciprocal sum can be written in the similar fashion
\begin{eqnarray}
E_{\rm int,k} = \frac{Z_{1}Z_{2}}{|\Gamma_1||\Gamma_2|}\int_{\Gamma_{1}}
\sum_{G_{1n}\in\RL_1^*,~G_{1n}\neq0}\frac{1}{|G_{1n}|^2}{\rm e}^{\frac{-|G_{1n}|^2}{4\eta^2}}{\rm e}^{-i G_{1n}\cdot ({\bf r}+{\bf t})}{\rm d}{\bf r} .
\end{eqnarray}
Since $G_{1n}\in\RL_1$ is a reciprocal lattice vector, which is orthogonal to ${\bf t}$, we have
\begin{eqnarray}
\int_{\Gamma_{1}} {\rm e}^{-i{G_{1n}\cdot ({\bf r}+{\bf t})}}{\rm d}{\bf r} = \int_{\Gamma_{1}} {\rm e}^{-i{G_{1n}\cdot {\bf r}}}{\rm d}{\bf r}=0 
\qquad{\rm for}~G_{1n}\in\RL_1,~G_{1n}\neq 0 .
\end{eqnarray}
This means that the averaged reciprocal sum $E_{\rm int,k}$ is $0$. 
Together with the Ewald sums within each lattice, we can compute the Ewald energy of the incommensurate system. 
Note that the above calculations can be easily extended to cases with multi atoms in the unit cell.


\small

\begin{thebibliography}{10}
	
	\bibitem{AAmodel}
	S.~Aubry and G.~Andr\'e.
	\newblock Analyticity breaking and anderson localization in incommensurate
	lattices.
	\newblock {\em Ann. Isr. Phys. Soc.}, 3:18, 1980.
	
	\bibitem{baake17}
	M.~Baake, D.~Damanik, and U.~Grimm.
	\newblock Aperiodic order and spectral properties.
	\newblock Snapshots of modern mathematics from Oberwolfach, 2015.
	
	\bibitem{blionv15}
	I.~V. Blinov.
	\newblock Periodic almost-{S}chr\"{o}dinger equation for quasicrystals.
	\newblock {\em Scientific Report}, 5:11492, 2015.
	
	\bibitem{britnell13}
	L.~Britnell, R.~M. Ribeiro, A.~Eckmann, R.~Jalil, B.~D. Belle, A.~Mishchenko,
	Y.~J. Kim, R.~V. Gorbachev, T.~Georgiou, S.~V. Morozov, A.~N. Grigorenko,
	A.~K. Geim, C.~Casiraghi, A.~H.~Castro Neto, and K.~S. Novoselov.
	\newblock Strong light-matter interactions in heterostructures of atomically
	thin films.
	\newblock {\em Science}, 340:1311--1314, 2013.
	
	\bibitem{cances17}
	E.~Canc\`{e}s, P.~Cazeaux, and M.~Luskin.
	\newblock Generalized {K}ubo formulas for the transport properties of
	incommensurate {2D} atomic heterostructures.
	\newblock {\em Journal of Mathematical Physics}, 58:06350, 2017.
	
	\bibitem{chen19}
	H.~Chen, A.~Zhou, and Y.~Zhou.
	\newblock in preparation, 2018.
	
	\bibitem{dingzhou09}
	J.~Ding and A.~Zhou.
	\newblock {\em Statistical Properties of Deterministic Systems}.
	\newblock Springer-Verlag, Berlin; Tsinghua University Press, 2009.
	
	\bibitem{ebnonnasir14}
	A.~Ebnonnasir, B.~Narayanan, S.~Kodambaka, and C.~V. Ciobanu.
	\newblock Tunable {MoS2} bandgap in {MoS2}-graphene heterostructures.
	\newblock {\em Applied Physics Letters}, 105:031603, 2014.
	
	\bibitem{geim13}
	A.~K. Geim and I.~V. Grigorieva.
	\newblock Van der waals heterostructures.
	\newblock {\em Nature}, 499:419, 2013.
	
	\bibitem{gross61}
	E.~P. Gross.
	\newblock Structure of a quantized vortex in boson systems.
	\newblock {\em Nuovo Cimento}, 20:454--477, 1961.
	
	\bibitem{helgaker00}
	T.~Helgaker, P.~Jorgensen, and J.~Olsen.
	\newblock {\em Molecular Electronic-Structure Theory}.
	\newblock Wiley, 2000.
	
	\bibitem{jiang14}
	K.~Jiang and P.~Zhang.
	\newblock Numerical methods for quasicrystals.
	\newblock {\em Journal of Computational Physics}, 256:428--440, 2014.
	
	\bibitem{koda16}
	D.~S. Koda, F.~Bechstedt, M.~Marques, and L.~K. Teles.
	\newblock Coincidence lattices of {2D} crystals: {H}eterostructure predictions
	and applications.
	\newblock {\em The Journal of Physical Chemistry C}, 120:10895--10908, 2016.
	
	\bibitem{komsa13}
	H.~P. Komsa and A.~V. Krasheninnikov.
	\newblock Electronic structures and optical properties of realistic transition
	metal dichalcogenide heterostructures from first principles.
	\newblock {\em Physical Review B}, 88:085318, 2013.
	
	\bibitem{Li2017i}
	X.~Li, X.~Li, and S.~Das Sarma.
	\newblock Mobility edges in one-dimensional bichromatic incommensurate
	potentials.
	\newblock {\em Phys. Rev. B}, 96:085119, 2017.
	
	\bibitem{Liu16}
	Y.~Liu, N.~O. Weiss, X.~Duan, H.~Cheng, Y.~Huang, and X.~Duna.
	\newblock Van der waals heterostructures and devices.
	\newblock {\em Nature}, 1:16042, 2016.
	
	\bibitem{loh15}
	G.~C. Loh and R.~Pandey.
	\newblock A graphene-boron nitride lateral heterostructure -- a
	first-principles study of its growth, electronic properties, and chemical
	topology.
	\newblock {\em Journal of Materials Chemistry C}, 3:5918--5932, 2015.
	
	\bibitem{lschenPRL}
	H.~P. L\"uschen, S.~Scherg, T.~Kohlert, M.~Schreiber, P.~Bordia, X.~Li, S.~Das
	Sarma, and I.~Bloch.
	\newblock Single-particle mobility edge in a one-dimensional quasiperiodic
	optical lattice.
	\newblock {\em Phys. Rev. Lett.}, 120:160404, 2018.
	
	\bibitem{martin04}
	R.~M. Martin.
	\newblock {\em Electronic Structure: Basic Theory and Practical Methods}.
	\newblock Cambridge University Press, 2004.
	
	\bibitem{massatt18}
	D.~Massatt, S.~Carr, M.~Luskin, and C.~Ortner.
	\newblock Incommensurate heterostructures in momentum space.
	\newblock {\em Multiscale Modeling \& Simulation}, 16:429--451, 2018.
	
	\bibitem{massatt17}
	D.~Massatt, M.~Luskin, and C.~Ortner.
	\newblock Electronic density of states for incommensurate layers.
	\newblock {\em Multiscale Modeling \& Simulation}, 15:476--499, 2017.
	
	\bibitem{Novo16}
	K.~S. Novoselov, A.~Mishchenko, A.~Carvalho, and A.~H.~Castro Neto.
	\newblock 2{D} materials and van der waals heterostructures.
	\newblock {\em Science}, 353, 2016.
	
	\bibitem{Settino2017}
	J.~Settino, N.~Lo Gullo, A.~Sindona, J.~Goold, and F.~Plastina.
	\newblock Signatures of the single-particle mobility edge in the ground-state
	properties of tonks-girardeau and noninteracting fermi gases in a bichromatic
	potential.
	\newblock {\em Phys. Rev. A}, 95:033605, 2017.
	
	\bibitem{Walter09}
	W.~Steurer and S.~Deloudi.
	\newblock {\em Crystallography of Quasicrystals: Concepts, Methods and
		Structures}.
	\newblock Springer-Verlag Berlin Heidelberg, 2009.
	
	\bibitem{Sun2015}
	M.~L. Sun, G.~Wang, N.~B. Li, and T.~Nakayama.
	\newblock Localization-delocalization transition in self-dual quasi-periodic
	lattices.
	\newblock {\em EPL (Europhysics Letters)}, 110:57003, 2015.
	
	\bibitem{terrones14}
	H.~Terrones and M.~Terrones.
	\newblock Bilayers of transition metal dichalcogenides: Different stackings and
	heterostructures.
	\newblock {\em Journal of Materials Research}, 29:373--382, 2014.
	
\end{thebibliography}
\end{document}